\documentclass{amsart}
\usepackage{amssymb}
\usepackage{color}

\usepackage{epsfig}
\usepackage{graphicx}
\usepackage{psfrag}

\usepackage{amsfonts}
\usepackage{amsmath}
\usepackage{url}

\newtheorem{theorem}{Theorem}
\newtheorem{sublemma}[theorem]{Sublemma}
\newtheorem{proposition}[theorem]{Proposition}
\newtheorem{property}[theorem]{Property}
\newtheorem{corollary}[theorem]{Corollary}
\newtheorem{lemma}[theorem]{Lemma}

\theoremstyle{remark}
\newtheorem{remark}[theorem]{Remark}
\newtheorem{observation}[theorem]{Observation}

\newtheorem{definition}[theorem]{Definition}

\newcommand{\matris}[4]{\left ( \!\! \begin{array}{cc}
#1 & #2 \\ #3 & #4 \end{array} \!\! \right )}
\newcommand{\1}{1\hspace{-2.5pt}{\rm l}}
\newcommand{\binomial}[2]{\left
(\!\!\begin{array}{c}#1\\#2\end{array}\!\!\right)}
\newcommand{\llbracket}{[\![}
\newcommand{\rrbracket}{]\!]}

\newcommand{\llabel}[1]{\label{#1}
}

\newcommand{\ledessin}[1]{\includegraphics{#1.eps}}
\begin{document}

\title[Punctured-torus groups]
{Triangulated cores of punctured-torus groups}
\author[F. Gu\'eritaud]{Fran\c{c}ois Gu\'eritaud}
\date{May 2006}

\begin{abstract}
We show that the interior of the convex core of a quasifuchsian punctured-torus group admits an ideal decomposition (usually an infinite triangulation) which is canonical in two different senses: in a combinatorial sense \emph{via} the pleating invariants, and in a geometric sense \emph{via} an Epstein-Penner convex hull construction in Minkowski space. The result extends to certain non-quasifuchsian punctured-torus groups, and in fact to all of them if a strong version of the Pleating Lamination Conjecture is true.
\end{abstract}

\maketitle

\section{Introduction}

Among Kleinian groups with infinite covolume, quasifuchsian groups, which are deformations of Fuchsian surface groups, are fundamental examples. It has been noted \cite{jorgensen} that such groups can be analyzed much more explicitly when the base surface is just a once-punctured torus --- the curve complex is then dual to a (locally finite) tree, which simplifies many issues. Punctured-torus groups, which still retain many features of general quasifuchsian groups, have thus become a favorite ``training ground'': Minsky's work \cite{minsky} on end invariants is an example. Similarly, Caroline Series \cite{series} was able to prove the Pleating Lamination Conjecture for quasifuchsian punctured-torus groups (general pleating laminations, in contrast, seem to pose tremendous technical challenges). In this paper, we prove that the convex core $V$ of any quasifuchsian punctured-torus group admits an ideal triangulation (or slightly more general decomposition) which relates ``as nicely as one could hope'' both to the intrinsic geometry of $V$, and to the combinatorics of the boundary pleatings of $V$. This answers several conjectures made in \cite{aswy1}. As a byproduct, we get enough control to re-prove the result of \cite{series}.

\subsection{Objects of study} \llabel{definelambda}
Let $S$ be the once-punctured torus endowed with its differential
structure. Let $\mathcal{C}$ be the set of homotopy classes of
simple closed curves in $S$. Choose in $\mathcal{C}$ a meridian $m$
and longitude $l$ whose intersection number is $1$. Fix an
identification $s$, called the \emph{slope}, between $\mathcal{C}$
and $\mathbb{P}^1{\mathbb Q}$, so that $s(m)=\infty$ and $s(l)=0$.
Let $(\alpha^+,\beta^+)$ and $(\alpha^-,\beta^-)$ be elements of
${\mathbb{R}^*}^2$ such that $\beta^+/\alpha^+$ and
$\beta^-/\alpha^-$ are distinct irrationals. Define the \emph{pleatings}
$\lambda^{\pm}:\mathcal{C} \rightarrow {\mathbb R}^{\pm}$ by
$$\lambda^+(\eta/\xi)=\left|\left|\begin{array}{cc} \beta^+& \eta \\
\alpha^+&\xi \end{array}\right|\right| ~~\text{ and }~~
\lambda^-(\eta/ \xi)=-\left|\left|\begin{array}{cc} \beta^-& \eta \\
\alpha^-& \xi \end{array}\right|\right|,$$ where the double bars denote
the absolute value of the determinant (and $\xi$, $\eta$ are coprime integers). Replacing a pair $(\alpha^{\pm}, \beta^{\pm})$ by its
negative does not change $\lambda^{\pm}$. Notice that
$\lambda^-<0<\lambda^+$ (as functions on $\mathcal{C}$).

\begin{theorem} \llabel{pureexistence}
There exists a cusped, non-complete hyperbolic
$3$-manifold $V$ homeomorphic to $S\times \mathbb{R}$ whose metric
completion $\overline{V}$ is homeomorphic to $S\times [0,1]=V\sqcup
S_{-\infty} \sqcup S_{+\infty}$, where $S_{-\infty}, S_{+\infty}$ are pleated surfaces whose pleating measures restrict to $\lambda^-$ and $\lambda^+$ (respectively) on $\mathcal{C}$. \end{theorem}
The bulk of this paper is devoted to producing such a $V$ (with some adjustments, the method also applies to rational pleatings $\lambda^{\pm}$). The puncture of $S$ is required to correspond to a cusp of $V$, so $\overline{V}$ will be the convex core of a quasifuchsian punctured-torus group $\Gamma \subset \text{Isom}^+(\mathbb{H}^3)$, i.e. the convex core of the manifold $\mathbb{H}^3/\Gamma$ (a punctured-torus group is a group freely generated by two elements with parabolic commutator). The interior $V$ of $\overline{V}$ is called the \emph{open convex core}.

Specifically, $V$ is constructed as an infinite union of
ideal hyperbolic tetrahedra $(\Delta_i)_{i\in \mathbb{Z}}$ glued
along their faces, and the $\lambda^{\pm}$ encode the gluing rule between $\Delta_i$ and $\Delta_{i+1}$ (see Section \ref{sectionstrategy}): so this ideal decomposition $\mathcal{D}$ of $V$ is canonical in a combinatorial
sense, with respect to the data $\lambda^{\pm}$. By \cite{series}, the group $\Gamma$ is determined uniquely up to conjugacy in $\text{Isom}(\mathbb{H}^3)$ by the $\lambda^{\pm}$. Our construction therefore provides a decomposition $\mathcal{D}$ of the open convex core $V$ of an \emph{arbitrary} quasifuchsian, non-Fuchsian punctured-torus group $\Gamma$. A good drawing of $V$ is Figure 3 of \cite{thurstonlimitedouble}.

\subsection{Geometric canonicity}
Akiyoshi and Sakuma \cite{comparing} generalized the Epstein-Penner convex
hull construction in Minkowski space $\mathbb{R}^{3+1}$ (see \cite{epsteinpenner}) to show that $V$ also admits a decomposition $\mathcal{D}^G$, canonical in a purely geometric sense, and related to the
Ford-Voronoi domain of $\mathbb{H}^3/\Gamma$. Roughly speaking, $\mathcal{D}^G$ is defined by considering the $\Gamma$-orbit $\mathcal{O}\subset \mathbb{R}^{3+1}$ of an isotropic vector representing the cusp, and projecting the cell decomposition of the boundary of the convex hull of $\mathcal{O}$ back to $\mathbb{H}^3/\Gamma$ (see Section \ref{sectioneph} for more detail). Our main theorem is

\begin{theorem} \llabel{ephgeneral}
If $V$ is the open convex core of a quasifuchsian once-punctured torus group,
the decompositions $\mathcal{D}$ and $\mathcal{D}^G$ of $V$ are the same.
\end{theorem}


It is reasonable to understand ending laminations of geometrically infinite surface groups as \emph{infinitely strong pleatings}, and to conjecture that the group is determined by its ending and/or pleating laminations.
Indeed, our method allows to construct (conjecturally unique) punctured-torus groups with arbitrary admissible ending and/or pleating laminations: the precise statement, with a full description of $\mathcal{D}^G$ (especially in the case of rational laminations), is Theorem \ref{ephhyper} in Section \ref{sectionextensions}.

\subsection{Context} The identity $\mathcal{D}=\mathcal{D}^G$ of Theorem \ref{ephgeneral} (and the existence of $\mathcal{D}$, as realized by positively oriented cells) was Conjecture 8.2 in \cite{aswy1}, also called the \emph{Elliptic-Parabolic-Hyperbolic (EPH) Conjecture}.
Akiyoshi subsequently established the identity in the case of two
infinite ends, in \cite{akiyoshi}. Near a finite end however, the ideal tetrahedra $\{\Delta_i\}_{i\in \mathbb{Z}}$ of $\mathcal{D}$ flatten at a very quick rate: the smallest angle of $\Delta_i$ typically goes to $0$ faster than any geometric sequence, as $i$ goes to $\pm \infty$. In a sense, the difficulty is to show that these angles nevertheless stay positive for the hyperbolic metric.

Under the hypotheses of Theorem \ref{ephgeneral}, for finitely many indices $i\in\mathbb{Z}$, the tetrahedron $\Delta_i$ of $\mathcal{D}^G$ comes from a \emph{spacelike} face in Minkowski space, and is therefore dual to a singular point (a vertex) of the Ford-Voronoi domain of $\mathbb{H}^3/\Gamma$. The latter domain is described in great detail in \cite{jorgensen, aswy1, aswy2} and \cite{aswy3}, relying on a \emph{geometric continuity} argument in the space of quasifuchsian groups (which is known to be connected). To study \emph{all} tetrahedra $\Delta_i$ at once, including those not seen in the Ford-Voronoi domain, the present paper takes a somewhat opposite approach: first describe geometric shapes for the tetrahedra of the candidate triangulation $\mathcal{D}$, then establish that the gluing of these tetrahedra defines (the open convex core of) a quasifuchsian group.

The construction of $\mathcal{D}$ (namely, of angles for the tetrahedra) will be fairly explicit: the solution will arise as the maximum of an explicit concave ``volume'' functional $\mathcal{V}$ over an explicit convex domain. The domain has infinite dimension, but there are explicit bounds on the contributions to $\mathcal{V}$ of the ``tail'' coordinates. This should allow for numerically efficient implementations.

The geometrically canonical decomposition $\mathcal{D}^G$ can be defined for arbitrary cusped manifolds, but is quite mysterious and hard to study in general. It seems to be unknown, for instance, whether $\mathcal{D}^G$ is always locally finite (see \cite{comparing}). For quasifuchsian punctured-torus groups however, Theorem \ref{ephgeneral} can be said to completely describe the combinatorics of $\mathcal{D}^G$. In fact, aside from coarsenings of $\mathcal{D}$, the author does not know of any ideal cell decomposition that is invariant under the hyperelliptic involution (a property which $\mathcal{D}^G$ must a priori enjoy).

An identity of the form $\mathcal{D}=\mathcal{D}^G$ can be established by the same methods in several closely related contexts: punctured-torus bundles (where the result is due to Lackenby \cite{lackenby}); complements of two-bridge links  (see the announcement in \cite{aswy3}) or of certain arborescent links --- see \cite{these} for a synthesis, but the present paper contains the key ideas. We will use several results (and notation) from \cite{mapomme}.

\subsection*{Acknowledgements}
My thanks go to Makoto Sakuma for having drawn my attention to the EPH conjecture. I am also deeply indebted to Francis Bonahon, Fr\'ed\'eric Paulin, Makoto Sakuma and David Futer for many discussions and insights.


\section{Strategy} \llabel{sectionstrategy}
\subsection{Setup} \llabel{subsectionsetup}
We return to the irrational pleating data $\lambda^{\pm}$.
In the hyperbolic plane ${\mathbb H}^2$ with boundary $\partial
{\mathbb H}^2={\mathbb P}^1{\mathbb R}$, consider the Farey
triangulation (the ideal triangle $01\infty$ iteratedly reflected in
its sides, see e.g. Section 2 of \cite{minsky}). The irrationals $\beta^+/\alpha^+$ and $\beta^-/\alpha^-$ belong to $\partial \mathbb{H}^2$, and the oriented line $\Lambda$ from $\beta^-/\alpha^-$ to $\beta^+/\alpha^+$ crosses infinitely many Farey edges $(e_i)_{i\in{\mathbb Z}}$ (the choice of $e_0$ is
arbitrary). To every pair of consecutive integers $i,i+1$ is
associated a letter, $R$ or $L$, according to whether $e_i$ and
$e_{i+1}$ share their Right or Left end, with respect to the
orientation of $\Lambda$ (we say that $\Lambda$ \emph{makes a Right} or
\emph{makes a Left} across the Farey triangle between $e_i$ and $e_{i+1}$).
We thus get a bi-infinite word $...RLLLRR...$ with infinitely many
$R$'s (resp. $L$'s) near either end.

Two rationals of $\mathbb{P}^1\mathbb{Q}$ are Farey neighbors exactly when
the corresponding elements (simple closed curves) of $\mathcal{C}$ have number of intersection $1$. Therefore, each Farey triangle $\tau$ defines an
\emph{ideal triangulation}
of the punctured torus $S$ in the following way. The
vertices of $\tau$ are the slopes of three curves of $\mathcal{C}$,
each parallel to a properly embedded line running from the puncture
to itself (so we may speak of the \emph{slopes} of such properly
embededd lines). These three lines are embedded disjointly and
separate $S$ into two \emph{ideal triangles}. Moreover, two
triangulations corresponding to Farey triangles which share an edge
differ by a \emph{diagonal move} (see Figure \ref{diagonalmove}).
Such a diagonal move must be seen as a (topological) ideal
tetrahedron in $S\times \mathbb{R}$ filling the space between two
(topological) pleated surfaces, pleated along the two ideal
triangulations.
\begin{figure}[h!] \centering
\ledessin{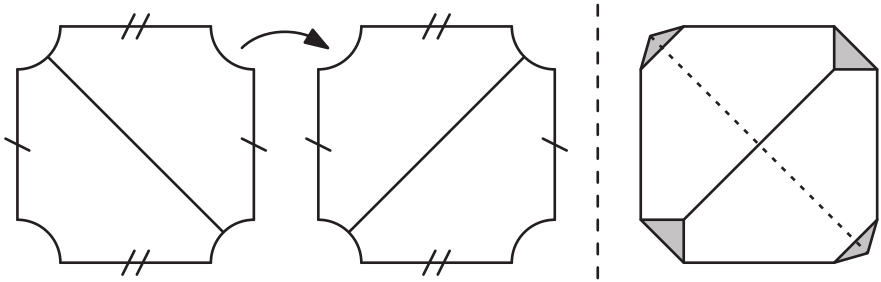}
\caption{\llabel{diagonalmove} Left: a diagonal move in $S$ (the puncture is in
the corners). Right: how to see it as an ideal tetrahedron, here with
truncated vertices (grey).} \end{figure}

Our strategy will be to consider triangulated surfaces (pleated
punctured tori) corresponding to the Farey triangles living between
$e_i$ and $e_{i+1}$ for $i$ ranging over ${\mathbb Z}$, interpolate
these surfaces with tetrahedra $\Delta_i$ corresponding to the
diagonal moves, and provide geometric parameters (dihedral angles)
for all these objects, using Rivin's Maximum Volume Principle (see
\cite{rivin}).

Throughout the paper, we will deal with an infinite family of
tetrahedra $(\Delta_i)_{i\in \mathbb{Z}}$, separated by pleated
once-punctured tori $S_i$. By an arbitrary choice, we resolve that
$S_i$ is the surface between $\Delta_i$ and $\Delta_{i-1}$ --- or
equivalently, that $\Delta_i$ is bounded by the surfaces $S_i$ and
$S_{i+1}$. However, the numbering of the tetrahedra should be seen
as the more essential one (see especially Definition
\ref{definitionhinge} below).

\subsection{Plan of the paper}
In Section \ref{sectionangles}, we describe the space of possible
dihedral angle assignments $x_i, y_i, z_i$ for the $\Delta_i$. In
Section \ref{sectionboundingthebending}, we encode $\lambda^{\pm}$
into constraints on the $x_i, y_i, z_i$. In Section
\ref{sectionhyperbolicvolume}, we carry out (constrained) volume
maximization. Important asymptotic features of the solution are analyzed in
Section \ref{sectionbehaviorofwi}. In Section
\ref{sectionthecusplink}, we describe the Euclidean triangulation of
the cusp. In Sections \ref{sectionintrinsicconvergence} and
\ref{sectionextrinsicconvergence}, we show that the pleated surfaces
$S_i$ converge in a strong enough sense, so that their limit as $i$
goes to $\pm \infty$ describes the (pleated) boundary of the metric
completion of $V=\bigcup_{i\in \mathbb{Z}} \Delta_i$. At that point,
we have constructed (the convex core of) a quasifuchsian group. The corresponding instance of Theorem \ref{ephgeneral} then follows from a computation, carried out in Section \ref{sectioneph}. In Section \ref{sectionextensions}, we provide a similar construction of punctured-torus groups with \emph{rational} pleating slopes $\beta^{\pm}/\alpha^{\pm}$ and/or with infinite ends, and re-prove that $(\lambda^+, \lambda^-)$ are continuous coordinates for the space of quasifuchsian groups (see \cite{series}).

\section{Dihedral angles} \llabel{sectionangles}

In this section we find positive dihedral angles for the ideal tetrahedra
$\Delta_i$, following Section 5 of \cite{mapomme}. More precisely, we describe the convex space $\Sigma$ of positive
angle configurations for the $\Delta_i$ such that:
\begin{itemize}
\item the three dihedral angles near each ideal vertex of $\Delta_i$ add up to $\pi$
(this is true in any ideal tetrahedron of ${\mathbb H}^3$);
\item the dihedral angles around any edge of $V=\bigcup_{i\in \mathbb{Z}}\Delta_i$
add up to $2\pi$ (this is necessary, though not sufficient, for a
hyperbolic structure at the edge);
\item the three pleating angles of each pleated punctured torus $S_i$ add up to $0$
(this is necessary, though not sufficient, to make the puncture $p$
of $S$ a cusp of $V$, i.e. make the loop around $p$ lift to a
parabolic isometry of $\mathbb{H}^3$).
\end{itemize}
(The first condition implies that opposite edges in $\Delta_i$ have
the same dihedral angle.) Later on we shall apply Rivin's Maximum Volume
Principle on a certain convex subset of (the closure of) $\Sigma$.

If the tetrahedron $\Delta_i$ realizes a diagonal exchange that
kills an edge $\varepsilon'$ and replaces it with $\varepsilon$,
denote by $\pi-w_i$ the interior dihedral angle of $\Delta_i$ at
$\varepsilon$ and $\varepsilon'$. Observe that the slope of
$\varepsilon$ (resp. $\varepsilon'$) is the rational located
opposite the Farey edge $e_i$ in the Farey diagram, on the side of
$\beta^+/\alpha^+$ (resp $\beta^-/\alpha^-$).

Thus, the pleating angles of the surface $S_i$ living between
$\Delta_{i-1}$ and $\Delta_i$ are \begin{equation}
\llabel{pleatingangles} w_{i-1}~~,~~-w_i~~\text{ and }~~
w_i-w_{i-1}~.\end{equation} Observe the sign convention: the pleated
punctured torus embedded in $S\times\mathbb{R}$ receives an
\emph{upward} transverse orientation from $\mathbb{R}$, and the
angles we consider are the dihedral angles \emph{above} the surface,
minus $\pi$. Thus, the ``new'' edge of $\Delta_{i-1}$, pointing
upward, accounts for a positive pleating $w_{i-1}$, while the
``old'' edge of $\Delta_i$, pointing downward, accounts for a
negative pleating $-w_i$. This is in accordance with the convention
$\lambda^-<0<\lambda^+$ of the Section \ref{definelambda}. One may write the three numbers (\ref{pleatingangles}) in the corners of the corresponding
Farey triangles (Figure \ref{rl2}, top).

In the tetrahedron $\Delta_i$, let $x_i$ (resp. $y_i$) be the
interior dihedral angle at the edge whose slope is given by the
right (resp. left) end of the Farey edge $e_i$. Let $z_i=\pi-w_i$ be
the third dihedral angle of $\Delta_i$. For instance, $2x_i$ (resp.
$2y_i$) is the difference between the numbers written just below and
just above the right (resp. left) end of $e_i$ in Figure \ref{rl2}
(the factor $2$ comes from the fact that the two edges of $\Delta_i$
with angle $x_i$ [resp. $y_i$] are identified). For notational
convenience, write $(w_{i-1},w_i,w_{i+1})=(a,b,c)$. By computing
differences between the pleating angles given in Figure \ref{rl2}
(bottom), we find the following formulae for $x_i, y_i, z_i$
(depending on the letters, $R$ or $L$, living just before and just
after the index $i$):

\begin{figure}[h!] \centering 
\psfrag{w}{$-w_i$}
\psfrag{w1}{$w_{i-1}$}
\psfrag{ww}{$w_i-w_{i-1}$}
\psfrag{e}{$e_i$}
\psfrag{e1}{$e_{i-1}$}
\psfrag{L}{$L$}
\psfrag{R}{$R$}
\psfrag{a}{$a$}
\psfrag{bb}{$-b$}
\psfrag{b}{$b$}
\psfrag{cc}{$-c$}
\psfrag{cb}{$c-b$}
\psfrag{ba}{$b-a$}
\psfrag{p}{$p$}
\psfrag{p1}{$p'$}
\psfrag{p2}{$p''$}
\ledessin{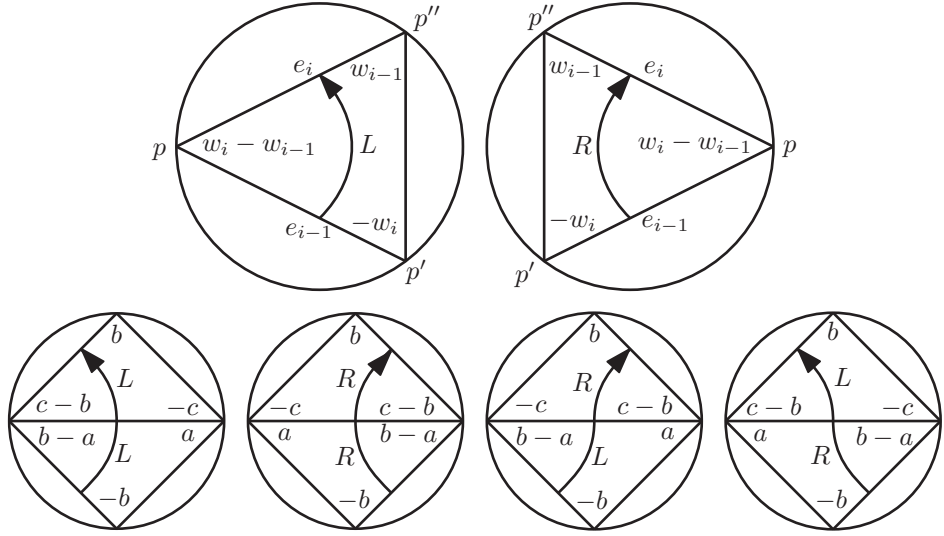}
\caption{Bottom: $e_i$ is the horizontal edge and
$(a,b,c)=(w_{i-1},w_i,w_{i+1})$.\llabel{rl2}}
\end{figure}

\begin{equation} \llabel{xiyi}
\begin{array}{c|c|c|c|c}
& L~~~~~L & R~~~~~R & L~~~~~R & R~~~~~L\\ \hline
x_i&\frac12 (a+c)&\frac12 (-a+2b-c)&\frac12 (a+b-c)&\frac12 (-a+b+c)\\
y_i&\frac12 (-a+2b-c)&\frac12 (a+c)&\frac12 (-a+b+c)&\frac12 (a+b-c)\\
z_i&\pi-b&\pi-b&\pi-b&\pi-b\\
\end{array}\end{equation}


The first of the three conditions defining $\Sigma$ can be checked immediately;
the other two are true by construction. From (\ref{xiyi}), the condition for
all angles to be positive is that:

\begin{equation} \llabel{positivity} \left \{ \begin{array}{l}
\text{ For all $i$ one has $0<w_i<\pi$.} \\
\text{ If $i$ separates identical letters (first two cases),
$2w_i>w_{i+1}+w_{i-1}$.} \\
\text{ If $i$ separates different letters (last two cases),
$|w_{i+1}-w_{i-1}|<w_i$.} \end{array} \right . \end{equation}
Denote by $\Sigma$ the non-empty, convex solution set of (\ref{positivity}). 

\section{Bounding the Bending} \llabel{sectionboundingthebending}

\subsection{A natural constraint on the pleating of $S_i$} \llabel{subsectionnaturalconstraint}
Next, we describe a certain convex subset of $\Sigma$. It is
obtained in the following way. Consider the pleated punctured torus,
$S_i$, lying between the tetrahedra $\Delta_i$ and $\Delta_{i-1}$.
Let $\epsilon_1, \epsilon_2,\epsilon_3$ be the edges of $S_i$ and
$\delta_1,\delta_2,\delta_3$ the corresponding pleating angles
(exterior dihedral angles, counted positively for salient edges as
in (\ref{pleatingangles}) above). Then, define the pleating measure
$\lambda_i:\mathcal{C}\rightarrow \mathbb{R}^+$ of $S_i$ by
\begin{equation} \lambda_i(\gamma)=\nu_1^{\gamma}\delta_1+\nu_2^{\gamma}\delta_2 +\nu_3^{\gamma}\delta_3, \llabel{definelambdai} \end{equation} where
$\nu_s^{\gamma}\in\mathbb{N}$ is the intersection number of $\gamma$ with the
simple closed curve $\varepsilon_s$ parallel to $\epsilon_s$. We
shall require that  $$\lambda^-(\varepsilon_s)<
\lambda_i(\varepsilon_s)< \lambda^+(\varepsilon_s)$$ for each
$s\in\{1,2,3\}$, and the same thing for every pleated punctured
torus $S_i$ in the $\mathbb{Z}$-family (\footnote{In fact a similar inequality will follow for all simple closed curves $\gamma$ of $\mathcal{C}$ (see Lemma \ref{bendingbounded} below), namely, $\lambda^-<\lambda_i<\lambda^+$ as functions on $\mathcal{C}$. Forcing this ``natural'' inequality
is the whole point of our constraint.}).

In other words, denote by $\eta/\xi$ the (rational) slope of
$\varepsilon_s$. Observe that
$\lambda_i(\varepsilon_s)=\delta_{s'}+\delta_{s''}=-\delta_s$ where
$\{s,s',s''\}=\{1,2,3\}$, because the slopes of $\varepsilon_1,
\varepsilon_2, \varepsilon_3$ are Farey neighbors. Therefore the
requirement is that, for every $\varepsilon_s$ as above,
\begin{equation} \llabel{encadrement} -\left|\left|
\begin{array}{cc} \beta^+&\eta\\
\alpha^+&\xi\end{array} \right|\right|< \delta_s
<\left|\left|\begin{array}{cc} \beta^-&\eta\\
\alpha^-&\xi\end{array} \right|\right|.\end{equation}

To express (\ref{encadrement}) in terms of the $w_i$, we need some
notation (the $\delta_s$ are determined by the $w_i$ via (\ref{pleatingangles}) above). For each Farey edge $e_i$, let $q_i^+$ (resp. $q_i^-$) be the rational
located opposite $e_i$, on the same side of $e_i$ as
$\beta^+/\alpha^+$ (resp. $\beta^-/\alpha^-$). For an arbitrary
rational $p=\eta/\xi$ (reduced form), introduce the slightly
abusive (\footnote{Abusive in that it depends on the ordered pair
$(\alpha,\beta)$ rather than just on the real
$\frac{\beta}{\alpha}$.}) notation \begin{equation} \llabel{wedge}
\frac{\beta}{\alpha} \wedge p := \left | \left | \begin{array}{cc} \beta&\eta\\
\alpha &\xi \end{array} \right | \right | \text{  (absolute value of
the determinant)}. \end{equation} Then, if $l,r$ are the rationals
living at the left and right ends of the Farey edge $e_i$, one has
\begin{equation} \llabel{signes} \textstyle{(\frac{\beta^+}{\alpha^+}
\wedge l ) + (\frac{\beta^+}{\alpha^+} \wedge r) =
\frac{\beta^+}{\alpha^+} \wedge q_i^-~~;~~
(\frac{\beta^-}{\alpha^-}\wedge l ) + (\frac{\beta^-}{\alpha^-}
\wedge r) = \frac{\beta^-}{\alpha^-} \wedge q_i^+}~.\end{equation}
Indeed, the $||\cdot||$-notation is invariant under $PSL_2({\mathbb
Z})$, acting on $({\mathbb H}^2, \partial {\mathbb H}^2)$ by
isometries and on ordered pairs $\pm(\alpha,\beta)$ as a matrix
group. Since $PSL_2({\mathbb Z})$ acts transitively on oriented
Farey edges, we are reduced to the case
$(l,r,q_i^+,q_i^-)=(\infty,0,1,-1)$ where $\beta^-/\alpha^-<0<
\beta^+/\alpha^+$, which is straightforward.

Let us now translate Equation (\ref{encadrement}) in terms of the
$w_i$. Let $e_{i-1}, e_i$ be two consecutive Farey edges; $p$ and
$p'$ are the ends of $e_{i-1}$; $p$ and $p''$ are the ends of $e_i$.
One has $q_{i-1}^+=p''$ and $q_i^-=p'$. In view of
(\ref{pleatingangles}) and Figure \ref{rl2} (top), Equation
(\ref{encadrement}) translates to $$\begin{array}{rcccl}
-(\frac{\beta^+}{\alpha^+} \wedge p'') &<&w_{i-1}&<&
\frac{\beta^-}{\alpha^-} \wedge p'' \\ -( \frac{\beta^+}{\alpha^+}
\wedge p' ) &<&-w_i&<& \frac{\beta^-}{\alpha^-} \wedge p' \\ -(
\frac{\beta^+}{\alpha^+} \wedge p ) &<&w_i-w_{i-1}&<&
\frac{\beta^-}{\alpha^-} \wedge p \\ \end{array}$$ Using the fact
that the $w_i$ are positive, this simplifies to (respectively)
$$\begin{array}{rcccl} &&w_{i-1}&<&\frac{\beta^-}{\alpha^-} \wedge q_{i-1}^+ \\
-(\frac{\beta^+}{\alpha^+} \wedge q_i^-) &<&-w_i&& \\
-(\frac{\beta^+}{\alpha^+} \wedge p) &<&w_i-w_{i-1}&<&
\frac{\beta^-}{\alpha^-} \wedge p \\ \end{array}$$

Finally, observe that $\frac{\beta^+}{\alpha^+} \wedge
p=(\frac{\beta^+}{\alpha^+} \wedge q_{i-1}^-) -
(\frac{\beta^+}{\alpha^+} \wedge q_i^-)$ while
$\frac{\beta^-}{\alpha^-} \wedge p=(\frac{\beta^-}{\alpha^-} \wedge
q_i^+) - (\frac{\beta^-}{\alpha^-} \wedge q_{i-1}^+)$, by Equation
(\ref{signes}). Therefore, if we introduce

\begin{equation}\llabel{phiphi} \phi_i^+:=\frac{\beta^+}{\alpha^+}
\wedge q_i^- ~~~\text{and }~~~
\phi_i^-:=\frac{\beta^-}{\alpha^-}\wedge q_i^+
\end{equation} then Equation (\ref{encadrement}) reduces to \begin{equation}
\left . \begin{array} {rcccl} w_i&<&\min\{\phi_i^+,\phi_i^-\} &&\\
\phi^-_{i-1}-\phi^-_i&<& w_{i-1}-w_i &<&\phi^+_{i-1}-\phi^+_i
\end{array} \right \}\forall i\in
\mathbb{Z}.\llabel{bounding}\end{equation}

\subsection{Study of $\phi^+$ and $\phi^-$} \llabel{subsectionstudyphi}

\begin{definition} \llabel{definitionhinge}
In \ref{subsectionsetup} we associated to each $i\in \mathbb{Z}$ a
Farey edge $e_i$ living between two letters of $\{R,L\}$. We call
$i$ a \emph{hinge index}, and $\Delta_i$ a \emph{hinge tetrahedron},
if the two letters are distinct ($RL$ or $LR$). \emph{Non-hinges}
correspond to $RR$ or $LL$.
\end{definition}

\begin{lemma} \llabel{studyphi} The following holds concerning the sequences $\phi^+, \phi^-:{\mathbb Z} \rightarrow {\mathbb R}_+^*$.
\begin{enumerate}
\item First, $\phi^-$ is strictly increasing and $\phi^+$
is strictly decreasing.

\item For all $i$ one has $1<\phi_{i+1}^-/\phi_i^-<2$ and
$1<\phi_{i-1}^+/\phi_i^+<2$.

\item If $i\in \mathbb{Z}$ is non-hinge then
$\phi_{i-1}^++\phi_{i+1}^+=2\phi_i^+$ and
$\phi_{i-1}^-+\phi_{i+1}^-=2\phi_i^-$.

\item If $i\in \mathbb{Z}$ is hinge then
$\phi_{i+1}^-=\phi_i^-+\phi_{i-1}^-$ and
$\phi_{i-1}^+=\phi_i^++\phi_{i+1}^+$.

\item The sequences $\phi^{\pm}$ are (weakly) convex, i.e. $2\phi^{\pm}_i\leq
\phi^{\pm}_{i-1}+\phi^{\pm}_{i+1}$.

\item If $i<j$ are consecutive hinge indices, then $1+\frac{j-i}{2} \leq
\{~\frac{\phi_i^+}{\phi_j^+}~,~ \frac{\phi_j^-}{\phi_i^-}~\} \leq
1+j-i$.

\item One has
$\underset{+\infty}{\lim}~\phi^+=\underset{-\infty}{\lim}~\phi^-=0$ and
$\underset{-\infty}{\lim}~\phi^+=\underset{+\infty}{\lim}~\phi^-=+\infty$.

\item One has $\underset{-\infty}{\lim}~(\phi^+_{i-1}-\phi^+_i) =
\underset{+\infty}{\lim}~(\phi^-_{i+1}-\phi^-_i)=+\infty$.
\end{enumerate}
\end{lemma}

\begin{remark} Points {\sf iii} and {\sf iv} are just equality cases of the
inequalities (\ref{positivity}), which encode positivity of the
angles $x_i, y_i, z_i$. \llabel{remarkphi}
\end{remark}

\begin{proof} We will deal only with $\phi^+$: the arguments for $\phi^-$ are
analogous. Let $e_{i-1}=pr_{i-1}$ and $e_i=pr_i$ be consecutive
Farey edges. We have $\phi_i^+=\frac{\beta^+}{\alpha^+} \wedge
r_{i-1}$ while $\phi_{i-1}^+=\frac{\beta^+}{\alpha^+}\wedge r_{i-1}
+ \frac{\beta^+}{\alpha^+}\wedge p$ by Equation (\ref{signes}),
hence {\sf i}.

For {\sf ii}, we need care only about the upper bound. Just observe
that $\phi_{i-1}^+=(\frac{\beta^+}{\alpha^+}\wedge r_i +
\frac{\beta^+}{\alpha^+}\wedge p)+\frac{\beta^+}{\alpha^+}\wedge p$
while $\phi_i^+=\frac{\beta^+}{\alpha^+}\wedge r_i +
\frac{\beta^+}{\alpha^+}\wedge p$.

For {\sf iii}, assume $e_{i+1}=pr_{i+1}$ so that
$\phi^+_k=\frac{\beta^+}{\alpha^+} \wedge p+\frac{\beta^+}{\alpha^+}
\wedge r_k$ for $|i-k|\leq 1$. For $k\in\{i,i+1\}$ the right hand
side is $\frac{\beta^+}{\alpha^+} \wedge r_{k-1}$ , so
$(\phi^+_{k-1}-\phi^+_k)=\frac{\beta^+}{\alpha^+} \wedge p$, which
is sufficient.

For {\sf iv}, assume $e_{i+1}=p'r_i$. In the notations of Formula
(\ref{phiphi}), we have $q_{i+1}^-=p$ and $q_i^-=r_{i-1}$, the ends
of $e_{i-1}$. This together with Equation (\ref{signes}) yields the
result. Point {\sf v} follows from {\sf iv} and {\sf ii} at hinges
indices, and from {\sf iii} at other indices. For Point {\sf vi},
observe that $\phi_i^+=\phi_j^++(j-i)\phi_{j+1}^+$ by {\sf iii-iv},
and conclude using {\sf ii}. Point {\sf vii} follows from {\sf vi},
and the presence of infinitely many hinge indices near either end.
Point {\sf viii} follows from {\sf vii}, {\sf v}, and {\sf iv}.
\end{proof}

\subsection{Behavior of the pleatings
$\lambda_i:\mathcal{C}\rightarrow\mathbb{R}$} \llabel{behavior}
For any real sequence $u\in \mathbb{R}^{\mathbb{Z}}$, define the
real sequence $\nabla u$ by $(\nabla u)_i=u_{i-1}-u_i$. Let us
summarize the conditions imposed on $w:{\mathbb Z}\rightarrow
{\mathbb R}_+^*$ (Eq. \ref{positivity} and \ref{bounding} above):
\begin{equation} \llabel{doubleve} \left \{ \begin{array}{rcccll}
0&<&w_i&<&\min\{\phi_i^+,\phi_i^-,\pi\}& \\
\nabla \phi^-_i&<&\nabla w_i&<&\nabla \phi^+_i
&\text{for all $i\in \mathbb{Z}$;} \\
&&|w_{i+1}-w_{i-1}|&<&w_i&\text{if $i$ is a hinge;} \\
&&w_{i+1}+w_{i-1}&<&2w_i&\text{otherwise.}
\end{array} \right .\end{equation}
It is a simple exercise to check that $w_i= \tanh\phi_i^+
\tanh\phi_i^-$, for instance, satisfies this system ($\tanh$ may be
replaced by any strictly concave monotonous function from ${\mathbb
R}_+$ to $[0,1)$ with the same $1$-jet at $0$). 
\begin{definition} If (\ref{doubleve}') denotes the system (\ref{doubleve}) in which all strong inequalities have been turned into weak ones, let $W\subset\mathbb{R^Z}$ be the solution space of (\ref{doubleve}'). \llabel{defineW}\end{definition}
Suppose $(w_i)_{i\in\mathbb{Z}}\in W$
and consider the corresponding pleating measures $\lambda_i:\mathcal{C}\rightarrow \mathbb{R}$ of $S_i$, the pleated punctured torus lying between the tetrahedra
$\Delta_i$ and $\Delta_{i-1}$.

\begin{lemma} \llabel{bendingbounded}
For any curve $\gamma \in \mathcal{C}$, the sequence
$(\lambda_i(\gamma))_{i\in\mathbb{Z}}$ is nondecreasing, with
$$\lambda^-(\gamma)~\leq~\underset{i\rightarrow -\infty}{\lim}\lambda_i(\gamma)
~\leq~\underset{i\rightarrow +\infty}{\lim}\lambda_i(\gamma)~\leq~\lambda^+(\gamma).$$
\end{lemma}
\begin{proof} First, observe that the pleating angles of $S_{i-1},S_i$ are
always of the form given in Figure \ref{monotonpleat}, where $x$
(resp. $y$) is positive, equal to the interior dihedral angle of the
tetrahedron at the horizontal (resp. vertical) edges of the square.
The closed curve $\gamma$ in the surface traverses the given
square a number of times, either vertically, or horizontally, or
diagonally (cutting off one of the four corners). The pleating along
$\gamma$ increases by $2y$ per vertical passage, $2x$ per horizontal
passage, and $0$ per diagonal passage, hence the monotonicity
statement. (This argument, or a variant of it, is also valid for
higher-genus surfaces and non-simple closed curves, as long as the
tetrahedron has positive angles.)
\begin{figure}[h!]
\centering 
\psfrag{xmy}{$-x-y$}
\psfrag{xpy}{$x+y$}
\psfrag{A}{$A$}
\psfrag{B}{$B$}
\psfrag{A2x}{$A-2x$}
\psfrag{B2y}{$B-2y$}
\psfrag{ga}{$\gamma$} 
\ledessin{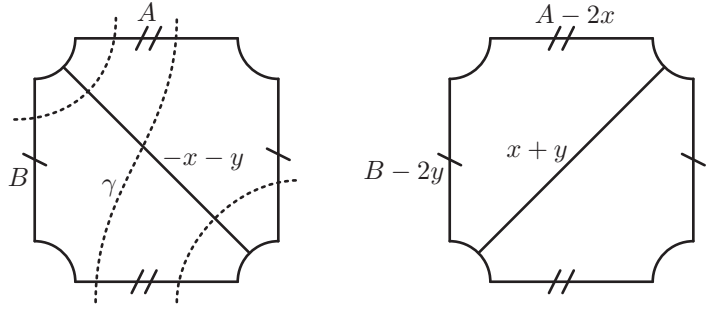}
\caption{The pleatings of
$S_{i-1}$ and $S_i$.\llabel{monotonpleat}}
\end{figure}

For the bounding, we will focus only on the positive side.
Consider the slope $s=\frac{\eta}{\xi}\in
\mathbb{P}^1\mathbb{Q}$ of $\gamma$ (a reduced fraction); recall the definition
$\lambda^+(s)=||^{\beta^+}_{\alpha^+} ~^{\eta}_{\xi}||$ from the
Introduction. Consider a large enough index $i$, such that the Farey edge $e_{i-1}$ separates $s$ from $\beta^+/\alpha^+$. Consider the points $p,p',p'' \in \mathbb{P}^1\mathbb{Q}$ such that $(e_{i-1},e_i)=(pp',pp'')$: the points $(s,p,\frac{\beta^+}{\alpha^+},p'',p',s)$ are cyclically arranged in $\mathbb{P}^1\mathbb{R}$. 
Observe that the $\wedge$-notation
(\ref{wedge}) applied at rationals is just the (geometric)
intersection number of the corresponding curves. Therefore, using
the angle information given in Figure \ref{rl2} (top), one has
\begin{equation}\llabel{determinant} \begin{array}{rcl} \lambda_i(\gamma)&=&(p\wedge s)(w_i-w_{i-1})-(p'\wedge s)w_i+(p''\wedge s)w_{i-1} \\ &=&  w_i (p\wedge s) + \nabla w_i (p' \wedge s) \hspace{25pt} \text{since $p''\wedge s=p\wedge s\,+\,p'\wedge s$} \\ &\leq &\phi^+_i(p\wedge s)+\nabla \phi^+_i(p'\wedge s)
\\ &=&\textstyle
{\left(\frac{\beta^+}{\alpha^+}\wedge p'\right)(p\wedge s)
+\left(\frac{\beta^+}{\alpha^+}\wedge p\right)(p'\wedge s)} \hspace{20pt}(*)
\end{array} \end{equation}
by definition (\ref{phiphi}) of $\phi^+$. The last quantity is $\frac{\beta^+}{\alpha^+}\wedge s$ (hence the upper bound): by $SL_2$-invariance of the $\wedge$-notation, it is enough to check this when $p=\infty$ and $p'=0$ --- in that case, $\frac{\beta^+}{\alpha^+}$ and $s=\frac{\eta}{\xi}$ have opposite signs, and $(*)$ does indeed
become $|\beta^+ \xi|+|\alpha^+ \eta|=\left | \left |^{\beta^+}_{\alpha^+} ~^{\eta}_{\xi} \right | \right |$.
\end{proof}

\section{Hyperbolic volume} \llabel{sectionhyperbolicvolume}

The product topology on $\mathbb{R^Z}$ induces a natural topology
on the space $W$ of Definition \ref{defineW}: clearly, $W$ is nonempty, convex, and compact.

If $(x,y,z)$ is a nonnegative triple such that $x+y+z=\pi$, let ${\mathcal
V}(x,y,z)$ be the hyperbolic volume of an ideal tetrahedron whose interior
dihedral angles are $x,y,z$. We wish to compute the total hyperbolic volume of
all tetrahedra when $w\in W$, i.e. $$\mathcal{V}(w):=\sum_{i\in\mathbb{Z}}
\mathcal{V}(x_i,y_i,z_i)$$ where $x_i,y_i,z_i$ are defined from the $w_i$
\emph{via} Table (\ref{xiyi}). This poses the problem of well-definedness --- the sum of the volumes might diverge. Let us estimate $\mathcal{V}$: a well-known explicit formula \cite{milnor} gives

\begin{equation} \begin{array}{rcl}
{\mathcal V}(x,y,z)&=&\displaystyle{-\int_0^x \log 2\sin - \int_0^y\log 2\sin - \int_0^z\log 2\sin} \\
&=&{\displaystyle \int_0^x \log \frac{\sin(\tau+y)}{\sin \tau} d\tau}
\hspace{15pt}\text{(as $\int_0^{\pi}\log 2 \sin = 0$)}\\ 
&\leq&\displaystyle{\int_0^x \log \frac{\tau+y}{\tau} d\tau} \hspace{20pt}
\text{(by concavity of $\sin$)}\\&=&\displaystyle{x\log \frac{x+y}{x} 
+y\log \frac{x+y}{y} } \\ &\leq& (x+y)\log 2 \hspace{35pt} \text{(concavity of
$\log$).} \end{array} \llabel{volumeformula} \end{equation}

\begin{lemma} \llabel{volumelemma} There exists a universal constant $K>0$ such
that the sum of the volumes of the tetrahedra $\Delta_j$ for $j\geq i$ (resp.
$j\leq i$) is at most $K\phi^+_i$ (resp. $K\phi^-_i$). \end{lemma}
\begin{proof} We will focus only on the $\phi^+$-statement. First, by the
computation above, the volume of the tetrahedron $\Delta_i$ is at
most $w_i \log 2 \leq \phi^+_i \log 2$ (see Table \ref{xiyi}). In
view of Lemma \ref{studyphi}-{\sf vi}, this implies that the total
volume of all \emph{hinge} tetrahedra beyond the index $i$ is at
most $3\phi^+_i\log 2$. For the same reason, it is sufficient to
prove
\begin{sublemma} \llabel{boundhops} There exists a
universal constant $L>0$ such that if $0$ and $N\in{\mathbb N}$ are two
consecutive hinge indices, then the sum of the volumes of the tetrahedra
$\Delta_1, \Delta_2, \dots, \Delta_{N-1}$ is at most $L \phi_0^+$.
\end{sublemma}

\begin{proof} In view of homogeneity in the
estimation (\ref{volumeformula}), it is sufficient to
assume $\phi_0^+=1$ and replace the volume with its estimate.
Therefore, let $(w_i)_{0\leq i \leq N}$ be a concave sequence in
$[0,1]$: following Table (\ref{xiyi}), we want to find a universal
upper bound $L$ (not depending on $N$) for $$
\sum_{i=1}^{N-1}\frac{w_{i+1}+w_{i-1}}{2}\log\frac{2w_i}{w_{i+1}+w_{i-1}}+
\frac{2w_i-w_{i+1}-w_{i-1}}{2}\log\frac{2w_i}{2w_i-w_{i+1}-w_{i-1}}.$$
If $A<B$ are positive integers, denote by $\Sigma_A^B$ the
restriction of the above sum to indices $A\leq i < B$. Observe that
the general term of $\Sigma_A^B$ is bounded by $2e^{-1}$, because
$\frac{\tau}{2}\log\frac{2}{\tau}\leq e^{-1}$ for all positive
$\tau$. To be more efficient, we bound the first half of the general
term by
$(\frac{w_{i+1}+w_{i-1}}{2})(\frac{2w_i}{w_{i+1}+w_{i-1}}-1)$, and
the second half, by concavity of $\log$. This produces
\begin{eqnarray*} \Sigma_A^B&\leq&
\sum_{i=A}^{B-1}\left (\frac{2w_i-w_{i+1}-w_{i-1}}{2} \times \frac{w_{i+1}+w_{i-1}}{w_{i+1}+w_{i-1}} \right ) \\
&& + \left ( \sum_{i=A}^{B-1} \frac{2w_i-w_{i+1}-w_{i-1}}{2}\right )
\log\frac{\sum_{i=A}^{B-1}w_i}{\sum_{i=A}^{B-1}
\frac{2w_i-w_{i+1}-w_{i-1}}{2}} \\ &=& \sigma \log \frac{e
\sum_{i=A}^{B-1}w_i}{\sigma}~,~~~~~\text{ where
}\sigma=\frac{w_A-w_{A-1} +w_{B-1}-w_B}{2}. \end{eqnarray*} Denote
by $M\in \llbracket 0, N \rrbracket$ a value of the index $i$ for
which $w_i$ is maximal. If $A<B\leq M$, we have $$0 \leq \sigma \leq
\frac{w_A-w_{A-1}}{2} \leq \frac{1}{2A} \leq 1$$ (the third
inequality follows from concavity of $w$ between $0$ and $A$). Also
notice that $f: \tau \mapsto \tau \log \frac{e}{\tau}$ is
nondecreasing on $[0,1]$. We shall apply these facts for $A=2^{k-1}$
and $B=\min\{2^k,M\}$: the previous bound on $\Sigma_A^B$ can be rewritten
\begin{eqnarray*} \Sigma_{A}^{B}&\leq&
\textstyle{ f(\sigma)+ \sigma \log\left ( \sum_{i=A}^{B-1} w_i \right )}\\
&\leq&\textstyle{f(2^{-k})+2^{-k}\log 2^{k-1} = 2^{-k}[1+(2k-1)\log
2].} \end{eqnarray*} The latter numbers (for $k$ ranging over
$\mathbb{N}^*$) add up to some universal $L'<+\infty$. After a
similar argument for the indices $M<i<N$, we can take
$L=2e^{-1}+2L'$.
\end{proof} Finally, we can take $K=3L+3\log 2$. Lemma \ref{volumelemma} is
proved. \end{proof}

\begin{corollary} The volume functional ${\mathcal V}:W\rightarrow {\mathbb
R}^+$ is well-defined, continuous, and concave.
\end{corollary}

\begin{proof} Well-definedness is the point of Lemma \ref{volumelemma}. Given
$\varepsilon>0$, only  finitely many indices $i$ satisfy
$\min\{\phi_i^+,\phi_i^-\}>\varepsilon/K$, and the others contribute at most
$2\varepsilon$ to the volume: hence continuity in the product topology.
Concavity follows from the concavity of the volume of \emph{one} tetrahedron
(parametrized by its angles): see e.g.
Proposition 8 of \cite{mapomme}. \end{proof}

Therefore, by compactness, there exists a sequence
$(w_i)_{i\in\mathbb{Z}}\in W$ which maximizes the hyperbolic volume
$\mathcal{V}$. From this point on, $w$ will denote that
maximizer.

\begin{proposition}
For each $j\in\mathbb{Z}$, if $x_j y_j z_j=0$ then $\max\{x_j,y_j,z_j\}=\pi$.
\llabel{tricot} \end{proposition}
\begin{proof}
Assume the tetrahedron $\Delta_j$ has exactly one vanishing angle, and aim for a contradiction. If $\mathcal{V}(\Delta^t)$ is the volume of a tetrahedron $\Delta^t$ having angles
$x^t,y^t,(\pi-x^t-y^t)$ with $(x^t, y^t)_{t\geq 0}$ smooth, $x^0=0<y^0<\pi$ and $dx^t/dt_{|t=0}>0$, then
$d\mathcal{V}(\Delta^t)/dt_{|t=0}=+\infty$ (by Formula (\ref{volumeformula}) above). 

Let $(w'_i)_{i\in\mathbb{Z}}$ be a sequence satisfying all (strict) inequalities of (\ref{doubleve}) and define $w^t:=w+t(w'-w)$ for $0\leq t \leq 1$. 
Denote by $\Delta^t_i$ the $i$-th tetrahedron determined via (\ref{xiyi}) by $(w^t_i)_{i\in\mathbb{Z}}$: then the angles of $\Delta^t_j$ satisfy the hypotheses above, so
$$\mathcal{V}(w^t)=\mathcal{V}(\Delta^t_j)+\left ( \sum_{i\neq j}\mathcal{V}(\Delta^t_i)\right )$$
has right derivative $+\infty$ at $t=0$ (the second summand is concave and continuous, so it has a well-defined right derivative in $\mathbb{R}\cup \{+\infty\}$ at $0$). Therefore, $\mathcal{V}$ was not maximal at $w$. \end{proof}


In Section \ref{sectionextensions} we will need the following consequence of Proposition \ref{tricot}:

\begin{proposition} \llabel{hingefragile}
If $j\in \mathbb{Z}$ and $x_j y_j z_j=0$, then $j$ is a hinge index and $w_j=0$. \qed \end{proposition}
This is Proposition 13 in \cite{mapomme}. But in fact much more is true:
\begin{proposition} \llabel{winterior}
All the (strict) inequalities of (\ref{doubleve}) are true at $w$.
\end{proposition}
\begin{proof}
If some inequality of (\ref{doubleve}) involving $\phi^+$
fails to be strict, it is easy to see
by induction that $w=\phi^+$ near $+\infty$, so all tetrahedra $\Delta_i$ (for $i$ large enough) have exactly one vanishing angle: $w$ was not maximal, by Proposition \ref{tricot}. Therefore all inequalities of (\ref{doubleve}) involving $\phi^{\pm}$ are strict. 
The arguments in \cite{mapomme} (especially
Lemma 16 and the argument of Section 9 there) can then be used
to show that \emph{all} inequalities (\ref{doubleve}) are strict at
$w$, so $w$ is a \emph{critical} point of the volume $\mathcal{V}$.
\end{proof}

Proposition \ref{winterior} implies that the holonomy representation
is trivial, i.e. the gluing of any finite number of consecutive
tetrahedra $\Delta_i$ defines a complete hyperbolic metric with
(polyhedral) boundary (the shapes of the $\Delta_i$ ``fit together
correctly'' around the edges): see \cite{rivin}, \cite{chanhodgson} or \cite{mapomme}. Therefore, the links of the vertices of the ideal tetrahedra
(Euclidean triangles) form a triangulation of the link of the
puncture: the latter is naturally endowed with a Euclidean structure
and its universal cover can be drawn in the plane (Figure
\ref{frise} --- more on the combinatorics of this triangulation in
Section \ref{sectionthecusplink}; see also \cite{mapomme}). Denote by $\Gamma$ the image of the induced holonomy representation $\pi_1(S)\rightarrow \text{Isom}^+(\mathbb{H}^3)$.
\begin{figure}[h!]
\centering 
\psfrag{L}{$\mathcal{L}$}
\ledessin{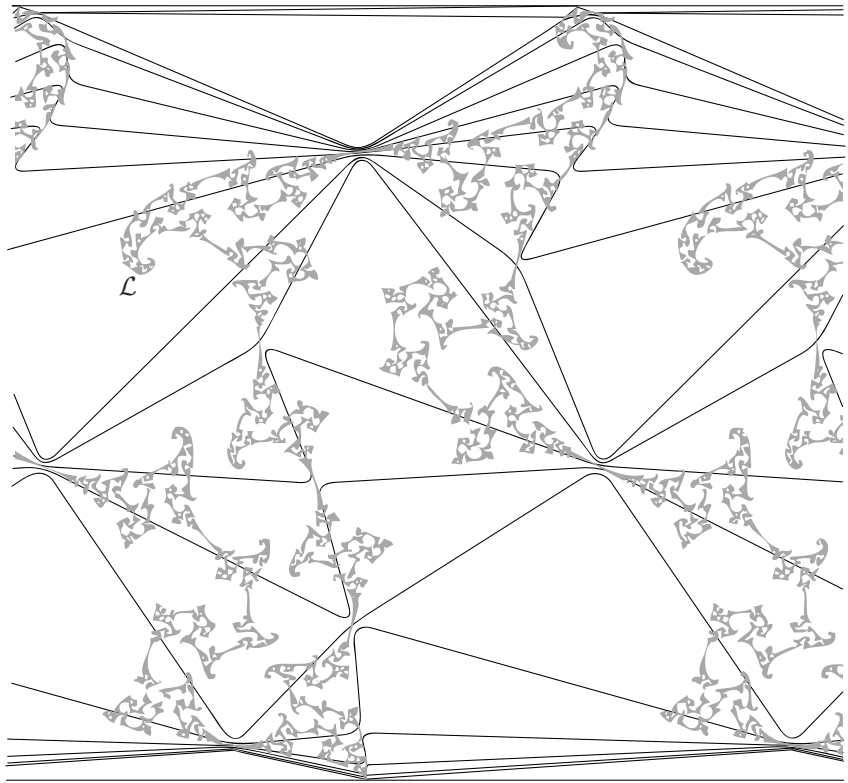}
\caption{The cusp triangulation is shown, in anticipation, against the limit set $\mathcal{L}\subset \mathbb{P}^1\mathbb{C}$ of the quasifuchsian group (a complicated Jordan curve). The picture extends periodically to the right and left. Each broken line is the puncture link of a pleated surface $S_i$; the vertices (whose design artificially sets the $S_i$ apart from each other) are all parabolic fixed points at which the limit set $\mathcal{L}$ becomes pinched. Infinitely many very flat triangles accumulate along the top and bottom horizontal lines. This picture was generated with Masaaki Wada's computer program Opti \cite{opti}. \llabel{frise}}
\end{figure}

\section{Behavior of $(w_i)$} \llabel{sectionbehaviorofwi}

\subsection{Properties of $(w_i)$}
We want to prove that the pleatings $\lambda_i$ (see Section
\ref{sectionboundingthebending}) of the pleated punctured tori $S_i$
converge (in the weak-* topology, i.e. on any test curve $\gamma \in
\mathcal{C}$) to the pleatings $\lambda^{\pm}$ near infinity. For
this we must study $(w_i)$, and especially show that the bounds
specified by $\phi^{\pm}$ in Equations (\ref{doubleve}) are almost
(but not quite) achieved. We use the $\nabla$-notation as in Section \ref{behavior}.

\begin{lemma} \llabel{hiatus} Recall that
$\nabla \phi^-<0<\nabla \phi^+$ (Lemma \ref{studyphi}). One has $$\max \left\{
\underset{\mathbb{Z}}{\overline{\lim}}\frac{w}{\phi^+}~,~
\underset{\mathbb{Z}}{\overline{\lim}}\frac{\nabla
w}{\nabla\phi^+}\right\} = \max
\left\{\underset{\mathbb{Z}}{\overline{\lim}}\frac{w}{\phi^-}~,~
\underset{\mathbb{Z}}{\overline{\lim}}\frac{-\nabla w}{-\nabla
\phi^-}\right\} = 1.$$ \end{lemma}

\begin{proof} We focus on the $\phi^+$-statement; the $\phi^-$-part is
analogous. Since $\frac{w}{\phi^+}<1$, assuming $\overline{\lim}
\frac{w}{\phi^+}<1$ implies $\sup \frac{w}{\phi^+}<1$, and the same
holds true for $\frac{\nabla w}{\nabla \phi^+}$ (see Equation
\ref{doubleve}). Therefore, suppose $\sup \frac{w}{\phi^+} \leq
1-\varepsilon$ and $\sup \frac{\nabla w}{\nabla \phi^+}\leq
1-\varepsilon$ for some $\varepsilon >0$, and aim at a
contradiction.

Recall the ordered pair $(\alpha^+, \beta^+)$ that helped define
$\phi^+$. For each $\mu>0$, define
$(\alpha^{\mu},\beta^{\mu}):=(\mu \alpha^+, \mu \beta^+)$. This
defines a new pleating function $\lambda^{\mu}=\mu \lambda^+ :
\mathcal{C} \rightarrow \mathbb{R}^+$, a new $\phi^{\mu}=\mu
\phi^+$ and a new domain $W^{\mu}$ by (\ref{doubleve}) (the numbers
$\alpha^-$ and $\beta^-$ are left unchanged). By definition, $W=W^1$
and $W^{\mu}\subset W^{\mu'}$ if and only if $\mu \leq \mu'$. Let
$(w^{\mu}_i)_{i\in \mathbb{Z}}$ be the maximizer of the volume
functional $\mathcal{V}$ on $W^{\mu}$. By assumption, we have
$w^1\in W^{1-\varepsilon}\subset W^1$, so $w^1=w^{1-\varepsilon}$.

Write $V(\mu)=\mathcal{V}(w^{\mu})$, so that
$V(1)=V(1-\varepsilon)$. It is straightforward to check that for any
$t\in [0,1]$, one has $(t\cdot w^{\mu}+(1-t)\cdot w^{\mu'}) \in
W^{t\mu+(1-t)\mu'}$. Since the volume of any tetrahedron is a
concave function of its angles, this is enough to imply that
$V:\mathbb{R}^+\rightarrow \mathbb{R}^+$ is (weakly) concave. By
inclusion, $V$ is also nondecreasing. In fact, $V$ is strictly
increasing (which will finish the proof by contradiction): to prove
this, since $V$ is concave, we just need to produce arbitrarily
large values of the volume $\mathcal{V}$ at points $v^{\mu} \in
W^\mu$, for large enough $\mu$. We may assume (up to a translation
of indices, see Lemma \ref{studyphi}, {\sf vii} -- {\sf viii}) that
$\phi^-_0>\pi$ and $-\nabla\phi^-_0>\pi$. Then, start by defining
$v^{\mu}_i=\min\{\phi^{\mu}_i,\phi^-_i,\pi\}$, so that
$v^{\mu}=\pi$ on $\llbracket 0,N+1\rrbracket$ for arbitrarily large
$N$. Without loss of generality, by just taking $\mu$ large enough,
we may further assume $|\nabla\phi^{\mu}|>\pi$ on $\llbracket
0,N+1\rrbracket$. Then, each time $i\in\llbracket 1,N\rrbracket$ is
a \emph{hinge} index, replace $v^{\mu}_i$ by $2\pi/3$: this is
allowed in $W^{\mu}$, by our assumptions on $\nabla\phi^-,
\nabla\phi^{\mu}$. The angles of the tetrahedron $\Delta_i$ are
then all in $\{\pi/2,\pi/3,\pi/6\}$ (see Table \ref{xiyi}). Since
there are infinitely many hinge indices near $+\infty$, the volume
can become arbitrarily large: QED. \end{proof}

\begin{corollary} \llabel{hiatus2} In fact, 
$\displaystyle{ \underset{i\rightarrow +\infty}{\overline{\lim}}
\min \left \{ \frac{w_i}{\phi^+_i} , \frac{\nabla w_i}{\nabla\phi^+_i} \right\}
= \underset{i\rightarrow -\infty}{\overline{\lim}}
\min \left \{ \frac{w_{i-1}}{\phi^-_{i-1}} , \frac{-\nabla
w_i}{-\nabla\phi^-_i} \right \} = 1 }$. \end{corollary}

\begin{proof} Again, we focus only on $\phi^+$. By Lemma \ref{hiatus}, there
exists a subsequence $(w_{\nu(i)})_{i\in\mathbb{N}}$ such that
$w_{\nu(i)}\sim\phi^+_{\nu(i)}$ or $\nabla w_{\nu(i)}\sim
\nabla\phi^+_{\nu(i)}$. Suppose the latter is the case. For an
arbitrary integer $i$, let $n$ be the smallest hinge index larger
than or equal to $\nu(i)$: observe that $\phi^+_{n+1}=\nabla\phi^+_n
=\nabla\phi^+_{\nu(i)}$ by Lemma \ref{studyphi}{\sf -iii-iv}, while
$w_{n+1}\geq \nabla w_n\geq \nabla w_{\nu(i)}$ by the positivity
conditions (\ref{doubleve}). Therefore, $\frac{w_{n+1}}
{\phi^+_{n+1}}\geq \frac{\nabla w_{\nu(i)}} {\nabla\phi^+_{\nu(i)}}$
so up to redefining $\nu$ we may assume simply
$w_{\nu(i)}\sim\phi^+_{\nu(i)}$.

Pick $\varepsilon>0$. Take $i$ such that \begin{equation}
\llabel{press}w_{\nu(i)}\geq(1-\varepsilon)\phi^+_{\nu(i)}.\end{equation}
Let $n$ be the smallest hinge index strictly larger than $\nu(i)$.
If $n=\nu(i)+1$ then  $$w_n\geq
w_{n-1}-\nabla\phi^+_n\geq(1-\varepsilon)
\phi^+_{n-1}-\phi^+_{n+1}=\phi^+_n-\varepsilon \phi^+_{n-1}\geq
(1-2\varepsilon)\phi^+_n;$$ $$\nabla w_n=w_{n-1}-w_n \geq (1-\varepsilon)\phi^+_{n-1}-\phi^+_n =\phi^+_{n+1}-\varepsilon\phi^+_{n-1}\geq
(1-3\varepsilon) \phi^+_{n+1}= (1-3\varepsilon)\nabla\phi^+_n$$
where Lemma \ref{studyphi} has been used several times. Therefore,
$\min\left \{\frac{w_n}{\phi^+_n}, \frac{\nabla w_n}{\nabla\phi^+_n}
\right\} \geq 1-3\varepsilon$.

If $n\geq \nu(i)+2$, we can find an index $k$ such that
$\frac{k-\nu(i)}{n-\nu(i)}\in[\frac12,\frac23]$. We will show that
$\min\left \{\frac{w_k}{\phi^+_k}, \frac{\nabla w_k}{\nabla\phi^+_k}
\right\} \geq 1-8\varepsilon$, which will finish the proof.

\noindent $\bullet$ By positivity of $w$ and concavity of $w$ between the points
$(\nu(i),k,n)$ one has $w_k\geq \frac13 w_{\nu(i)}$. Therefore,
$$\phi^+_k-w_k\leq \phi^+_{\nu(i)}-w_{\nu(i)}\leq
\left (\frac{1}{1-\varepsilon}-1\right )w_{\nu(i)} \leq
\frac{3\varepsilon}{1-\varepsilon}w_k$$ (here the first inequality
holds because $(\phi^+-w)$ is decreasing by Condition
(\ref{doubleve}), and the second follows from the assumption
(\ref{press}) above). Hence,
$$\frac{w_k}{\phi_k^+}\geq\frac{1-\varepsilon}{1+2\varepsilon}\geq
1-3\varepsilon.$$

\noindent $\bullet$ Observe that $$(1-\varepsilon)\phi^+_{\nu(i)}\leq w_{\nu(i)}\leq w_k+(k-\nu(i))\nabla w_k \leq \phi^+_k+(k-\nu(i))\nabla w_k$$ where
the second inequality follows from concavity of $w$ between the
points $(\nu(i),k-1,k)$. It follows that $$(k-\nu(i))\nabla w_k\geq
\phi^+_{\nu(i)}-\phi^+_k-\varepsilon \phi^+_{\nu(i)}
=(k-\nu(i))\nabla\phi^+_k-\varepsilon \phi^+_{\nu(i)}$$ hence $$
\nabla w_k \geq
\nabla\phi^+_k-\frac{\varepsilon}{k-\nu(i)}\phi^+_{\nu(i)} \geq
\nabla\phi^+_k-\frac{2\varepsilon}{n-\nu(i)}\phi^+_{\nu(i)}.$$ Since
Lemma \ref{studyphi}-{\sf vi} gives us $\phi^+_{\nu(i)}\leq
2(n-\nu(i))\phi^+_n$, this yields $$\nabla w_k \geq
\nabla\phi^+_k-4\varepsilon \phi^+_n \geq
\nabla\phi^+_k-8\varepsilon \phi^+_{n+1}=(1-8\varepsilon)\nabla
\phi^+_k.$$
\end{proof}

\begin{corollary} Recall the pleating $\lambda_i$ of the pleated surface $S_i$.
For any simple closed curve $\gamma\in \mathcal{C}$ we have $\underset{i\rightarrow +\infty}{\lim}\lambda_i(\gamma)=\lambda^+(\gamma)$ and
$\underset{i\rightarrow -\infty}{\lim}\lambda_i(\gamma)=\lambda^-(\gamma)$. 
\end{corollary}
\begin{proof}
By Corollary \ref{hiatus2}, the member ratio in the inequality (\ref{determinant}) can be made arbitrarily close to $1$.
\end{proof}

\section{The cusp link} \llabel{sectionthecusplink}

We now aim to investigate the behavior of the pleated surfaces $S_i$
as $i$ goes to $\pm \infty$ --- or, more precisely, to find two
limiting pleated surfaces $S_{\pm\infty}$ with pleatings
$\lambda^{\pm}$ such that $\overline{V}=\bigcup_{i\in\mathbb{Z}}
\Delta_i \sqcup S_{+\infty} \sqcup S_{-\infty}$ is metrically
complete with locally convex boundary. The difficult part is to
prove that the intrinsic moduli of the $S_i$ converge in
Teichm\"uller space. This question will be addressed in the next
section. In this section, we just describe the cusp link, introduce
notation and prove a few inequalities.

As always, we will mainly work near $+\infty$. Define
$V=\bigcup_{i\in\mathbb{Z}}\Delta_i$. We begin by orienting all the
edges of $V$ in a way that will be consistent with the pleating data
$(\alpha^{\pm},\beta^{\pm})$. Namely, denote by $\mathcal{E}
\subset {\mathbb P}^1{\mathbb Q}$ the collection of all the
endpoints of all the Farey edges $(e_i)_{i\in {\mathbb Z}}$. For
$s\in \mathcal{E}$, let $Q_s$ be the edge of $V$ of slope $s$.
Recall that the punctured torus is defined as $({\mathbb
R}^2\smallsetminus {\mathbb Z}^2)/{\mathbb Z}^2$. Orienting $Q_s$ is
therefore equivalent to orienting the line $F_s$ of slope $s$ in
$\mathbb{R}^2$. We decide that the positive half of $F_s$ should be
on the same side of the line ${\mathbb R} (\alpha^-, \beta^-)$ as
$(\alpha^+, \beta^+)$, and orient $Q_s$ accordingly.

Let us now describe the link of the puncture, or cusp triangulation.
Each tetrahedron $\Delta_i$ contributes four similar Euclidean
triangles at infinity to the link of the puncture, corresponding to
the four ideal vertices of $\Delta_i$. The bases of these four
triangles form a closed curve which is a broken line of four
segments, and the triangles point alternatively up and down from
this broken line (see Figure \ref{fourtriangles}). The two upward
(resp. downward) pointing triangles have the same Euclidean size,
an effect of the hyperelliptic involution (rotation of $180^{\circ}$
around the puncture) which acts isometrically on $V$ (reversing all
edge orientations) and as a horizontal translation on the cusp link.

\begin{figure}[h!] \centering 
\psfrag{Q0}{$Q_0$}
\psfrag{Q1}{$Q_1$}
\psfrag{Qinf}{$Q_{\infty}$}
\psfrag{inf}{$(\!\bigstar\!)$}
\psfrag{tof}{: toward $(\!\bigstar\!)$}
\psfrag{frf}{: from $(\!\bigstar\!)$}
\psfrag{x}{$x$}
\psfrag{y}{$y$}
\psfrag{z}{$z$}
\psfrag{bap}{$\displaystyle{\frac{\beta^+}{\alpha^+}}$}
\ledessin{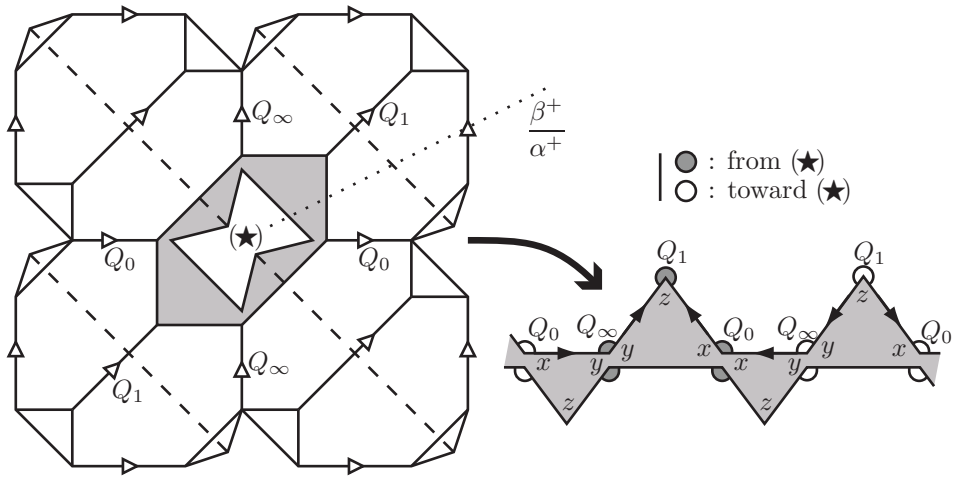}
\caption{In the left panel, $(\!\bigstar\!)$ marks the cusp. In the right panel, the cusp is at infinity. \llabel{fourtriangles}}
\end{figure}

\begin{definition} \llabel{definebases}
If the loop around the cusp has Euclidean length $4$, let $b_i$
(resp. $b'_i$) be the length of the base of a downward-pointing
(resp. upward-pointing) triangle contributed by the tetrahedron
$\Delta_i$, and define $\sigma_i=b'_i/b_i$ (see e.g. Figure \ref{sigmalimit}).
\end{definition}

Figure \ref{fourtriangles} also shows some additional information,
assuming that $e_i=0\infty$ and $\alpha^+, \beta^+>0$ (hence
$\beta^-/\alpha^-<0$). Namely, the pleated surfaces $S_{i+1}$ and
$S_i$ (above and below $\Delta_i$) are pleated along $Q_0,
Q_{\infty}, Q_1$ and $Q_0, Q_{\infty}, Q_{-1}$ respectively, and the
orientations of the $Q_s$ are as shown in the left panel of Figure
\ref{fourtriangles} (edge $Q_{-1}$ has no determined orientation).
The orientations of the lines to/from the puncture are also shown in
the right panel, at the vertices, with the help of a color code.
Moreover, each segment $\epsilon$ of the upper broken line in the
right panel corresponds to an arc about a vertex $v$ of a face $f$
of $S_{i+1}$ in the left panel, so $\epsilon$ receives the
orientation of the edge of $f$ opposite $v$. If $\tau$ is one of the
upward-pointing triangles drawn in the plane $\mathbb{C}$ (right
panel), consider the tetrahedron $\Delta$ whose vertices are
$\infty$ and those of $\tau$: all edges of $\Delta$, except one,
receive orientations from the construction above, and $\Delta$ is
isometric (respecting these orientations) to $\Delta_i$. Finally,
notice the labels in the $3$ corners of each triangle in the right
panel: the corner of the free vertex is labeled $z$ (the angle there
being $z_i$); the other two corners are labeled $x$ and $y$
accordingly. The labels $x$-$y$-$z$ appear clockwise in each triangle.

The contribution of the tetrahedron $\Delta_{i-1}$ to the link at
infinity is also a union of four triangles bounded by two broken
lines. Moreover, the upper broken line from $\Delta_{i-1}$ is the
lower broken line from $\Delta_i$, and the orientations of the lines
to/from infinity must agree. Inspection shows that there are only
two possibilities, corresponding to whether the letter living on the
surface $S_i$ is $R$ or $L$: for $R$, the $x$-corners of the two
$4$-triangle families live near a common vertex $C$; for $L$, the same
is true of the $y$-corners (Figure \ref{eighttriangles}).
\begin{figure}[h!] \centering 
\psfrag{tof}{: to $\infty$}
\psfrag{frf}{: from $\infty$}
\psfrag{A}{$A$}
\psfrag{B}{$B$}
\psfrag{C}{$C$}
\psfrag{b}{$b$}
\psfrag{c}{$c$}
\psfrag{R}{$(R)$}
\psfrag{L}{$(L)$}
\ledessin{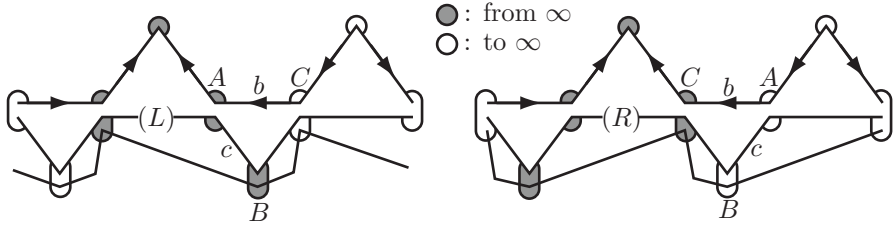}
\caption{Left and Right transitions. Compare with the right panel of Figure \ref{fourtriangles}. \llabel{eighttriangles}} \end{figure}

\begin{definition}
In a downward-pointing triangle defined by $\Delta_i$, the edge
lengths are $b_i$ (the basis from Def. \ref{definebases}), $b_{i-1}$, and a third number which we call $c_i$.
\end{definition}

\begin{property} For $i$ large enough, $(b_i)$ is increasing and $(b'_i)$ is
decreasing. \end{property}
\begin{proof} Let $\tau_i$ (resp. $\tau'_i$) be a downward-pointing
(resp. upward-pointing) triangle defined by the tetrahedron
$\Delta_i$. For large enough $i$ we have $z_i=\pi-w_i>\pi-\phi_i^+
>\pi/2$, so $b_i$ is the longest edge of $\tau_i$, and $b'_i$ the longest edge of $\tau'_i$. Since $b_{i-1}$
is an edge of $\tau_i$ and $b'_{i+1}$ is an edge of $\tau'_i$
(Figure \ref{eighttriangles} or \ref{sigmalimit}), the conclusion follows.
\end{proof}
\begin{property}We have $\lim_{i\rightarrow +\infty}\sigma_i=0$.
\end{property}
\begin{proof}
We already know that $(\sigma_i)=(b'_i/b_i)$ is ultimately
decreasing. It is therefore enough to show that
$\sigma_{i+1}/\sigma_{i-1}\leq 1/2$ for large enough hinge indices
$i$. Consider Figure \ref{sigmalimit},
where angles labeled $z$ are obtuse (as a rule, we shade a cusp triangle contributed by $\Delta_i$ whenever $i$ is a hinge index). Check that
$$\frac{\sigma_{i+1}}{\sigma_{i-1}}= \frac{b'_{i+1}}{b'_{i-1}}
\frac{b_{i-1}}{b_{i+1}} < \frac{b'_{i+1}}{b'_i} \frac{b_{i-1}}{b_i}=
\frac{\sin x_i \sin y_i}{\sin^2 z_i}\leq \frac{\sin^2(w_i/2)}{\sin^2
w_i}<\frac12$$ (the last two inequalities follow from an easy study
of $\sin$, using $x_i+y_i=\pi-z_i=w_i<\pi/2$). As an immediate
consequence, we find $\lim_{+\infty}b'_i=0$ and
$\lim_{+\infty}b_i=2$. Similarly, $\lim_{-\infty}b'_i=2$ and
$\lim_{-\infty}b_i=0$. (Compare with Figure \ref{frise}).
\begin{figure}[h!] \centering
\psfrag{xi}{$x_i$}
\psfrag{yi}{$y_i$}
\psfrag{zi}{$z_i$}
\psfrag{zip}{$z_{i+1}$}
\psfrag{zim}{$z_{i-1}$}
\psfrag{bi}{$b_i$}
\psfrag{bip}{$b_{i+1}$}
\psfrag{bim}{$b_{i-1}$}
\psfrag{Bi}{$b'_i$}
\psfrag{Bip}{$b'_{i+1}$}
\psfrag{Bim}{$b'_{i-1}$}
\ledessin{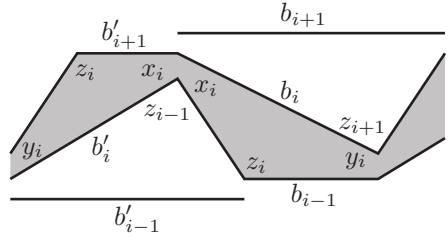}
\caption{The index $i$ is a hinge; the nature of $i\pm 1$ is undetermined. \llabel{sigmalimit}} \end{figure}
\end{proof}
\begin{definition}
Let $J\subset {\mathbb Z}$ be the set of all integers $j$ such that
$j-1$ is a hinge index. \llabel{definitionJ}
\end{definition}
\begin{proposition} \llabel{ecrase}
If $j<l$ are two large enough consecutive elements of $J$, and $k$
is not in $J$, then $$\frac{c_k}{c_{k-1}}= \sigma_{k-1}~~~
\text{and} ~~~ \frac{c_l}{c_j} \leq \sigma_{j-2}.$$
\end{proposition}
\begin{proof}
Since $k-1$ is not a hinge index, $b_{k-1}$ shares the same end with
$b_k$ and with $b_{k-2}$
(see Figure \ref{ecrasefigure}, left) so
$c_k/c_{k-1}=b'_{k-1}/b_{k-1}=\sigma_{k-1}$.

Since $c_l$ is always an edge of $\tau'_{l-1}$, we have $c_l\leq
b'_{l-1}$, hence
\begin{equation} \frac{c_l}{c_j}\leq\frac{b'_{l-1}}{c_j}
\leq\frac{b'_j}{c_j}
=\frac{c_{j-1}}{b_{j-2}}\leq\frac{b'_{j-2}}{b_{j-2}}
=\sigma_{j-2}\end{equation} where the equality in the middle just
translates the similarity property of the ``hinge'' triangles
$\tau_{j-1},\tau'_{j-1}$
(shaded in Figure \ref{ecrasefigure}, right).
\end{proof}
\begin{figure}[h!] \centering
\psfrag{cj}{$c_j$}
\psfrag{bj}{$b_j$}
\psfrag{Bj}{$b'_j$}
\psfrag{cjm}{$c_{j-1}$}
\psfrag{bjm}{$b_{j-1}$}
\psfrag{Bjm}{$b'_{j-1}$}
\psfrag{bjM}{$b_{j-2}$}
\psfrag{BjM}{$b'_{j-2}$}
\psfrag{ck}{$c_k$}
\psfrag{ckm}{$c_{k-1}$}
\psfrag{bkm}{$b_{k-1}$}
\psfrag{Bkm}{$b'_{k-1}$}
\ledessin{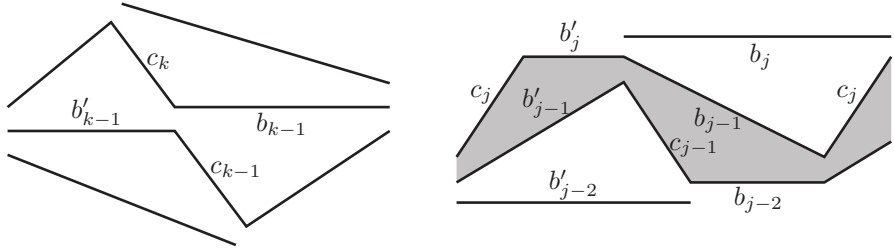}
\caption{Left: $k-1$ is not a hinge. Right: $j-1$ is a hinge. \llabel{ecrasefigure}}
\end{figure}

\section{Intrinsic convergence of the surfaces $S_i$} \llabel{sectionintrinsicconvergence}

\subsection{Thickness of the tetrahedra}

Consider a tetrahedron $\Delta_i$ bounded by the pleated surfaces
$S_i$ and $S_{i+1}$. Let $Q$ (resp. $Q'$) be the pleating edge of
$S_i$ (resp. $S_{i+1}$) not lying in $S_i \cap S_{i+1}$. Let $s_i$
be the shortest segment between $Q$ and $Q'$, across $\Delta_i$.

\begin{definition}
Recall the orientations on the edges of $V=\bigcup_{i\in
\mathbb{Z}}\Delta_i$. Let $\ell_i$ be the complex length of the
hyperbolic loxodromy along $s_i$ sending $Q$ to $Q'$, respecting the
orientations of $Q, Q'$ (with $-\pi< {\rm Im} \, \ell_i \leq
\pi$).
\end{definition}

\begin{proposition} The series $(\ell_i)_{i\in {\mathbb Z}}$ is
absolutely convergent. \llabel{helicoptere} \end{proposition}

\begin{proof} Consider a downward-pointing triangle $\tau$ contributed by
$\Delta_i$. Label the vertices of $\tau$ by $A,B,C$ in such a way
that $AC=b_i$, $BC=b_{i-1}$ and $AB=c_i$ (see Figure
\ref{eighttriangles}). Let $\gamma_i$ be the hyperbolic loxodromy of
complex length $\ell_i$ along the common perpendicular to $B\infty$
and $AC$, sending $B,\infty$ to $A,C$ (in that order).

\begin{sublemma} Let $\ell=\ell_i=\rho+\theta\sqrt{-1}$ be the complex length of
$\gamma_i$, with $\theta\in]-\pi,\pi]$. Then
$\max\{\rho,|\theta|\}\leq\pi\sqrt{c_i/b_i}$.
\llabel{thickness}\end{sublemma} \begin{proof} Up to a plane
similarity, we may assume $A=1$ and $C=-1$. Also, for convenience, relabel the edges of $ABC$ by $a,b,c$. Let $L$ be the fixed
line of $\gamma_i$. The hyperbolic isometry defined by $z\mapsto
f(z)= \frac{2B+1-z}{z+1}$ exchanges the oriented lines $AC$ and
$B\infty$, so it reverses the orientation of $L$ around the center
of the tetrahedron $ABC\infty$. Therefore, $\gamma_i$ is given by
$\gamma_i(z)= 2B-f(z)$. If $M$ is a matrix of $\gamma_i$, one
has $\displaystyle{\frac{ \text{tr}^2 M}{4\det
M}=\frac{\cosh \ell +1}{2}}$. Use
$M=\matris{2B+1}{-1}{1}{1}$ to find $B=\cosh \ell$. Compute
\begin{eqnarray*} \frac{a\pm c}{b}&=&\frac{|B+1|\pm
|B-1|}{2}=\left|\cosh^2 \frac{\ell}{2}\right|\pm \left|\sinh^2
\frac{\ell}{2}\right|
\\&=&\textstyle{ \cosh \frac{\ell}{2}\cosh
\frac{\overline{\ell}}{2}\pm\sinh \frac{\ell}{2} \sinh
\frac{\overline{\ell}}{2}= \cosh\frac{\ell\pm\overline{\ell}}{2}}.
\end{eqnarray*} Use $\cosh(i\theta)=\cos \theta$ to get $\cosh(\rho)=
\frac{a+c}{b}$ and $\cos(\theta)=\frac{a-c}{b}$. The estimates
$\text{Argcosh}(y)\leq 2 \sqrt{\frac{y-1}{2}}$ and $\text{Arccos}
(y) \leq \pi \sqrt{\frac{1-y}{2}}$, since $\frac12 \left |\frac{a\pm
c}{b}-1 \right |\leq
\frac{c}{b}$, finally yield $\rho\leq 2\sqrt{c/b}$ and $|\theta|\leq
\pi\sqrt{c/b}$, hence Sublemma \ref{thickness}.\end{proof}

Since $(b_i)$ goes to $2$, Proposition
\ref{helicoptere} will follow if the $c_i$ go to $0$ fast enough
near $+\infty$ (with a similar argument near $-\infty$). Such fast
decay is given by Proposition \ref{ecrase}: using the fact that
$(\sigma_i)$ goes to $0$, we can bound the $c_i$ by decreasing
geometric sequences on intervals of the form $\llbracket j, l-1
\rrbracket$ where $j<l$ are consecutive elements of the set $J$.
\end{proof}

\subsection{Teichm\"uller charts}

Let $\mathcal{T}_S$ be the Teichm\"uller space of the punctured
torus: for each $i\in {\mathbb Z}$ we denote by $\mu(S_i)\in
\mathcal{T}_S$ the intrinsic modulus of the marked surface $S_i$. In
order to prove that the $\mu(S_i)$ converge in $\mathcal{T}_S$, let
us first introduce appropriate charts for $\mathcal{T}_S$.

Consider a topological ideal triangulation $E$ of the punctured torus $S$,
with labeled edges $\epsilon_1,\epsilon_2,\epsilon_3$. Then $E$ defines an isomorphism
$$h_E: {\mathbb P}^2 {\mathbb R}_+^* =
(\mathbb{R}_+^*)^{\{\epsilon_1,\epsilon_2,\epsilon_3\}}/\mathbb{R}_+^*~\widetilde{\longrightarrow}~
\mathcal{T}_S.$$ Namely, given a hyperbolic metric $g$ on $S$, in
order to compute $h_E^{-1}(g)$, straighten $E$ to an ideal
triangulation for $g$ and return the positive projective triple of
Euclidean lengths defined (in the link of the puncture) by the
sectors opposite $\epsilon_1,\epsilon_2,\epsilon_3$. We consider the $h_E$ as charts of $\mathcal{T}_S$. We endow $\mathbb{P}^2\mathbb{R}_+^*$ with the
distance $d$ given by
$$d([a:b:c],[a':b':c']):=\min_{\lambda>0}\max\left \{ \left | \log
\frac{\lambda a}{a'}\right |, \left | \log \frac{\lambda
b}{b'}\right |, \left | \log \frac{\lambda c}{c'}\right | \right
\}.$$ Also define $\overline{h}_E:\mathbb{P}^2{\mathbb
C}^*\rightarrow \mathcal{T}_S$ by
$\overline{h}_E([a:b:c]):=h_E([\,|a|:|b|:|c|\,]).$


In particular, if the pleated punctured torus $S_{i+1}$, pleated
along the ideal triangulation $E_{i+1}$, gives rise in the cusp link
to a broken (oriented) line whose segments have complex coordinates
$(a,b,c)$, then $$\mu(S_{i+1})=\overline{h}_{E_{i+1}}([a:b:c]).$$

\begin{figure}[h!] \centering
\psfrag{x}{$a$}
\psfrag{y}{$b$}
\psfrag{z}{$c$}
\psfrag{xy}{$a+b$}
\psfrag{Y}{$\displaystyle{\frac{b}{a+b}c}$}
\psfrag{X}{$\displaystyle{\frac{a}{a+b}c}$}
\psfrag{Si}{$S_i$}
\psfrag{Sp}{$S_{i+1}$}
\ledessin{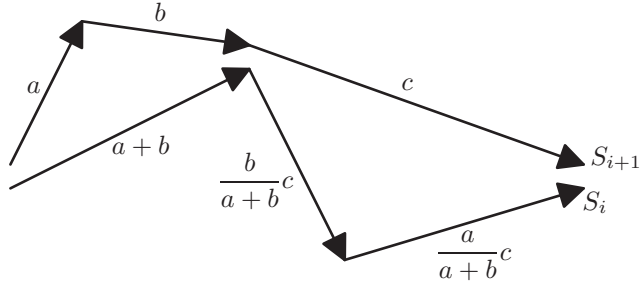}
\caption{The complex projective triples associated to $S_{i+1}$ and $S_i$.\llabel{chartmap}}
\end{figure}

Moreover, Figure \ref{chartmap} shows \emph{two} broken lines corresponding
respectively to $S_{i+1}$ and the previous pleated surface $S_i$.
Since the triangles (links of vertices of $\Delta_i$) are similar,
the complex coordinates of the segments forming the lower broken
line are functions of $(a,b,c)$, as shown. Therefore, if we define
the substitution formula
$$\Psi([a:b:c]):=\left (\left [
a+b~:~\frac{b}{a+b}c~:~\frac{a}{a+b}c \right ] \right )$$ (which is
a birational isomorphism from $\mathbb{P}^2\mathbb{C}$ to itself),
then $$\mu(S_{i})=\overline{h}_{E_{i}} (\Psi([a:b:c])).$$ Finally,
by viewing the diagonal exchange between $E_i$ and $E_{i+1}$ as a
\emph{flat} tetrahedron (as in Figure \ref{diagonalmove}), we see
that
$$\mu(S_{i+1})=h_{E_i}(\Psi([\,|a|:|b|:|c|\,])).$$ In other words, the
restriction $\psi$ of $\Psi$ to $\mathbb{P}^2\mathbb{R}_+^*$ is the
chart map $h_{E_{i+1}}\rightarrow h_{E_i}$, i.e.
$$\begin{array}{rcl}
\mathbb{P}^2\mathbb{R}_+^* & \overset{\psi}{\longleftarrow} & \mathbb{P}^2\mathbb{R}_+^* \\
h_{E_i} \searrow && \swarrow h_{E_{i+1}} \\
&\mathcal{T}_S& \end{array}$$ commutes.

\begin{property} \llabel{skid} We have $\displaystyle{d\left ( h_{E_i}^{-1}\mu(S_i), h_{E_i}^{-1}\mu(S_{i+1}) \right ) = \log \frac{|a|+|b|}{|a+b|}}$.
\end{property}
\begin{proof} The right member is
$$d\left(\left[ |a+b|:\left |\frac{b}{a+b}c \right |: \left |\frac{a}{a+b}c \right |\, \right],
\left [|a|+|b|:\frac{|b|}{|a|+|b|}|c|:\frac{|a|}{|a|+|b|}|c|\right
]\right ).$$ \end{proof} It is easy to see that
$\psi:\mathbb{P}^2\mathbb{R}_+^*\rightarrow
\mathbb{P}^2\mathbb{R}_+^*$ is $3$-Lipschitz for $d$. In fact,

\begin{proposition} There exists $K>0$ such that the $n$-th
iterate $\psi^n$ is $Kn$-bilipschitz for all $n>0$. \llabel{magique}
\end{proposition} \begin{proof} Since
$\psi([a:b:c])=\left [\frac{(a+b)^2}{c}:b:a\right ]$ and
$\psi^{-1}([a:b:c])=\left [c:b:\frac{(c+b)^2}{a}\right ]$, it is
enough to check the Lipschitz statement. Set
$A:=\sqrt{a}~;~C:=\sqrt{c}~$;
$$P_n:=\frac{A^{n+1}}{C^n}+\sum_{i,j\in\mathbb{Z}}\binomial{j}{i-1}
\binomial{n-i}{j-i} A^{2i-n-1} C^{n-2j}~~\text{ for all $n\geq
-1$}$$ (the sum is really on $0<i\leq j \leq n$), so that
$P_0=A~,~P_{-1}=C~,~P_1= \frac{A^2+1}{C}$. We claim that
$\psi^n[a:1:c]=[P_n^2:1:P_{n-1}^2]$. The claim is seen by induction
on $n$: the only difficult thing is the induction step
$P_{n+1}^2=\frac{(P_n^2+1)^2}{P_{n-1}^2}$. First, it is
straightforward to check that
$P_{n+1}+P_{n-1}=(\frac{A}{C}+\frac{C}{A}+\frac{1}{AC})P_n$ (using
the Pascal relation twice). Hence, for all $n\geq 1$, one has
\begin{eqnarray*} &&
(P_{n+1}P_{n-1}-P_n^2)-(P_n P_{n-2}-P_{n-1}^2)\\&=&
P_{n-1}(P_{n+1}+P_{n-1})-P_n(P_n+P_{n-2})\\&=&(P_{n-1}P_n-P_n P_{n-1})
\textstyle{ \left (\frac{A}{C}+\frac{C}{A}+\frac{1}{AC}\right )}=0.
\end{eqnarray*} Therefore $P_{n+1}P_{n-1}-P_{n}^2=P_1 P_{-1}-P_0^2=1$, which
proves the induction step. Since $P_n$ is a Laurent polynomial in
$A,C$ with partial degrees of order $n$ and positive coefficients,
we see that $\log P_n$ is $Ln$-bilipschitz in $\log a, \log c$ for
some universal $L$. The Proposition follows.
\end{proof}
The proof of Proposition \ref{magique} may seem extremely \emph{ad
hoc} and unsatisfactory. However, Proposition \ref{magique} is a
special case of a more general phenomenon for \emph{Markoff maps}
(in the sense of \cite{bowditch}), which we describe in
\cite{abeilles}.

\subsection{Convergence of the moduli}

\begin{proposition} The moduli $(\mu(S_i))_{i\rightarrow \pm \infty}$ converge in Teichm\"uller space $\mathcal{T}_S$.
\llabel{surfaceconvergence} \end{proposition}

\begin{proof} To prove the $+\infty$-statement, we will fix a large
enough index $i$, and prove that the series $$\eta_j:= d\left
(h_{E_i}^{-1}\mu(S_j), h_{E_i}^{-1}\mu(S_{j+1}) \right ),$$ defined
for $j> i$, has finite sum.

For $j>i$, consider a downward-pointing triangle $\tau$ contributed
by $\Delta_j$, with its edge lengths $b_j$, $b_{j-1}$ and $c_j$. The
angles of $\tau$ at the ends of $b_j$ are $x_j$ and $y_j$. By
Property \ref{skid}, we have $d\left
(h_{E_j}^{-1}\mu(S_j),h_{E_j}^{-1}\mu(S_{j+1})\right )= \log
\frac{c_j+b_{j-1}}{b_j}$. Compute \begin{eqnarray*}\log
\frac{c_j+b_{j-1}}{b_j}&\leq& \frac{c_j+b_{j-1}-b_j}{b_j}=\frac{\sin
x_j+\sin y_j-\sin z_j}{\sin z_j}=\frac{2\sin \frac{x_j}{2} \sin
\frac{y_j}{2}}{\cos\frac{x_j+y_j}{2}}\\
&\leq& \frac{\sin x_j \sin y_j}{\sin^2 z_j} \sin^2 z_j = \frac{c_j
b_{j-1}}{b_j^2} \sin^2 w_j \leq \frac{c_j}{b_j} {\phi^+_j}^2
\end{eqnarray*} (the inequality at the start of the second line
holds for large enough $j$ because $x_j+y_j=2w_j\rightarrow 0$).
Define $\delta_j=\frac{c_j}{b_j}{\phi_j^+}^2$, and let $M_j$ denote
the best bilipschitz constant for the chart map $h_{E_j}\rightarrow
h_{E_i}$. Observe that $\eta_j\leq M_j \delta_j$.

Let $j<l$ be two consecutive elements of $J$ (see Definition
\ref{definitionJ}), and $k\notin J$ an integer. We shall bound the
$M_n\delta_n$ by geometric sequences on intervals of the form
$\llbracket j,l-1\rrbracket$, using Proposition \ref{ecrase} as in
the proof of Proposition \ref{helicoptere}. Since $\psi$ is
$3$-bilipschitz, one clearly has $M_k\leq 3 M_{k-1}$. By Proposition
\ref{ecrase},
$$\frac{M_k \delta_k}{M_{k-1}\delta_{k-1}} \leq 3
\frac{\delta_k}{\delta_{k-1}} \leq 3 \frac{c_k}{c_{k-1}} =
3\sigma_{k-1}$$ because $(\phi^+_k)$ and $(1/b_k)$ are decreasing.
The right member goes to $0$ for large $k$.

By Sublemma \ref{magique}, $M_{l}\leq (l-j)L\cdot M_j$ for some
universal $L$, and by Lemma \ref{studyphi} ({\sf i-iii-iv}),
$\phi^+_{j}=\phi^+_{l-1}+(l-j-1)\phi^+_l\geq(l-j)\phi^+_{l}$. Using
Proposition \ref{ecrase}, it follows that
$$\frac{M_l\delta_l}{M_j\delta_j}\leq (l-j)L \cdot \frac{c_l}{c_j}
{\left ( \frac{\phi^+_l}{\phi^+_j} \right )}^2 \leq
\frac{L\sigma_{j-2}}{l-j}\leq L\sigma_{j-2}.$$ The right member goes
to $0$ as $j$ goes to infinity, hence Proposition
\ref{surfaceconvergence}.
\end{proof}

\section{Extrinsic convergence of the surfaces $S_i$}
\llabel{sectionextrinsicconvergence}

\subsection{Pleated surfaces} \llabel{definepleating} Propositions
\ref{helicoptere} and \ref{surfaceconvergence}, together with Lemma
\ref{bendingbounded}, are the key ingredients to prove that the
metric completion of $V=\bigcup_{i\in \mathbb{Z}}\Delta_i$ has two
boundary components which are pleated punctured tori with pleating
measure $\lambda^{\pm}$. A pleated surface is by definition (up to taking a universal cover) a map
$\varphi:\mathbb{H}^2\rightarrow \mathbb{H}^3$ which sends
rectifiable arcs to rectifiable arcs of the same length, such that
through each point $p$ of $\mathbb{H}^2$ runs an open segment $s_p$
on which $\varphi$ is totally geodesic. It is known (see \cite{lms},
5.1.4) that the direction of $s_p$ is unique if and only if $p$
belongs to a certain \emph{geodesic lamination} $\Lambda$ (closed
union of disjoint geodesics), and that $\varphi$ is totally geodesic
away from $\Lambda$.

To wrestle with pleated surfaces, we will use the fact that if
$\varphi$ is a locally convex immersion, and $\Lambda$ has zero
Lebesgue measure, then $\Lambda$ comes with a transverse (pleating)
measure $\nu_{\Lambda}$. More precisely, $\nu_{\Lambda}$ can be
defined on any segment $s$ transverse to $\Lambda$ in the following
way (see Sections 7 to 9 of \cite{bonahon}). Immerse $\varphi(\mathbb{H}^2)$ into the Poincar\'e upper half-space model. Each component of $\mathbb{H}^2\smallsetminus \Lambda$ crossed by $s$ can be extended to a subset $A$ of
$\mathbb{H}^2$ bounded by only one or two lines of $\Lambda$ crossed by
$s$. Endow $s$ with a transverse orientation. The boundary component
of $A$ in $\partial \mathbb{H}^3=\mathbb{C}\cup\{\infty\}$ on the
positive side of the transverse orientation defines a circle arc
$c_A$, of angle $\theta_A \in (-2\pi,2\pi)$ (we may assume
$\infty\notin c_A$). The closure of the union of all the $c_A$ forms
a rectifiable curve $c$, of length $\sum_A\text{\sf length}(c_A)$.
Then $c$ has a well-defined \emph{regular curvature} $RC(c)$,
defined as the absolutely convergent sum $\sum_A\theta_A$. But $c$
also has a \emph{total curvature} $TC(c)$, defined (up to an
appropriate multiple of $2\pi$) as the difference between the
arguments of the initial and final tangent vectors to $c$. (The appropriate multiple of $2\pi$ can be determined by closing off the \emph{embedded} arc $c$ with a broken line, and requiring that the resulting Jordan curve have total curvature $2\pi$). Then, $\nu_{\Lambda}(s)$ is defined as the \emph{singular curvature} $SC(c)=TC(c)-RC(c)$.

Conversely, if $\varphi:\mathbb{H}^2\rightarrow \mathbb{H}^3$ is a
pleated immersion and $SC(c)$ is well-defined and non-negative for
all transverse segments $s$, then $\varphi$ is locally convex with
pleating measure $\nu_{\Lambda}$ as above. We refer to \cite{bonahon} for greater detail.

\subsection{Setup} Consider the marked once-punctured torus
$S_{+\infty}$, endowed with the hyperbolic metric $\lim_{i\rightarrow
+\infty}\mu(S_i)$. There exists a unique compact geodesic lamination
$\Lambda_c^+$ of slope $\beta^+/\alpha^+$ on $S_{+\infty}$. Exactly
two lines $\ell,\ell'$ issued from the puncture of $S_{+\infty}$ fail
to meet $\Lambda_c^+$: define $\Lambda^+=\Lambda_c^+\sqcup \ell
\sqcup \ell'$. Then $S_{+\infty}$ is the disjoint union of
$\Lambda^+$ and the interiors of two ideal triangles $A,A'$. The
union $A\cup A'\cup \ell \cup \ell'$ is a punctured ideal bigon.

Recall the ideal tetrahedra $\Delta_i$, the space $V=\bigcup_{i\in \mathbb{Z}} \Delta_i$ and the
hyperelliptic involution $h:V\rightarrow V$ which reverses all edge
orientations. The slopes of the (oriented) pleating lines of the
surface $S_i$ (between $\Delta_{i-1}$ and $\Delta_i$) are
elements of ${\mathbb P}^1{\mathbb Q}$ projecting to $0,1,\infty$ in
${\mathbb P}^1({\mathbb Z}/2)$: accordingly, we call these pleating lines $l^0_i, l^1_i, l^{\infty}_i$. For $*\in\{0,1,\infty\}$, denote by
$\omega^*_i$ the unique point of $l^*_i$ fixed under the
hyperelliptic involution (the $\omega^*_i$ are called Weierstrass
points). Let $s^*_i$ be the segment from $\omega^*_i$ to
$\omega^*_{i+1}$ (across the tetrahedraon $\Delta_i$ if
$\omega^*_i\neq\omega^*_{i+1}$): each $s^*_i$ is contained in a
(pointwise) fixed line $\Omega^*$ of the hyperelliptic involution.

Fix the value of the superscript $*$ (soon we shall omit it). Fix a
point $\omega$ in $\mathbb{H}^2$ and an oriented line $l$  through
$\omega$. For each (oriented) surface $S_i$, consider the oriented, marked
universal covering $\pi_i:(\mathbb{H}^2,\omega,l)\rightarrow
(S_i,\omega^*_i,l^*_i)=(S_i,\omega_i,l_i)$. Endow a universal
covering $\widetilde{V}$ of $V$ with lifts $\widetilde{\omega_i}$ of
the $\omega_i$ connected by lifts of the segments $s_i$, and fix a
developing map $\Phi:\widetilde{V}\rightarrow \mathbb{H}^3$. There
is a unique map $h_i$ such that
\begin{equation} \llabel{comdiagram} \begin{array}{ccccc}
(\mathbb{H}^2,\omega,l)&\overset{h_i}{\longrightarrow}&(\widetilde{V},\widetilde{\omega_i})
& \overset{\Phi}{\longrightarrow}&\mathbb{H}^3 \\ \pi_i \downarrow
~~ && \downarrow &&\\ S_i & \longrightarrow & V && \end{array}
\end{equation}
commutes. We will prove that the (developing) pleated immersions
\begin{equation}\varphi_i=\Phi \circ h_i : \mathbb{H}^2 \rightarrow
\mathbb{H}^3 \llabel{immersion}\end{equation} converge as pleated
maps.

\subsection{Convergence of the $\varphi_i$} \llabel{hausdorff}
By Proposition \ref{helicoptere}, it is already clear that the
restriction of $\varphi_i$ to $l$ converges to a totally geodesic
embedding of the line $l$ into $\mathbb{H}^3$ (the convergence is
uniform on all compacts of $l$). By Ascoli's theorem, since the
$\varphi_i$ are $1$-Lipschitz, there exists an increasing sequence
$\nu$ such that the $\varphi_{\nu(i)}$ converge to a certain map
$\varphi_{+\infty}$, uniformly on all compact sets of $\mathbb{H}^2$. 
In Section \ref{pleatingphiinf} below, $\varphi_{+\infty}$ is shown to be independent of the subsequence $\nu$: in anticipation, we now abusively write $\varphi_i$ instead of $\varphi_{\nu(i)}$.

We shall work in the projective tangent bundle
$\mathcal{E}=\mathbb{P}T\mathbb{H}^2$, a circle bundle over
$\mathbb{H}^2$ in which geodesic laminations naturally live as
closed sets. Length and angle measurements define a (canonical)
complete Riemannian metric on $\mathcal{E}$. For $K\subset
\mathbb{H}^2$ compact and $A,B\subset \mathcal{E}$ closed, let
$K^{\mathcal E}\subset \mathcal{E}$ be the preimage of $K$ under the
natural projection $\mathcal{E}\rightarrow \mathbb{H}^2$, and define
$$d_K(A,B):=\inf \left \{ \delta>0 ~|~
A\cap K^{\mathcal E} \subset B+\delta ~,~ B\cap K^{\mathcal E}
\subset A+\delta \right \}$$ where $X+\delta$ denotes the set of
points within $\delta$ of $X$. Then $\inf(1,d_K)$ is a pseudometric,
and the set of closed subsets of $\mathcal{E}$ is compact for the
Hausdorff metric
$$d_H=\sum_{n>0}2^{-n}\inf(1,d_{K_n})$$ where the $K_n$ are concentric balls
of radius $n$.

\medskip
Observe that $S_{+\infty}$ also has Weierstrass points
$\omega^*_{+\infty}$, belonging to leaves $l^*_{+\infty}$ of the
lamination $\Lambda^+$. Fixing the value of $*$ as before, denote by
$\pi_{+\infty}:(\mathbb{H}^2,\omega,l) \rightarrow (S_{+\infty},
\omega^*_{+\infty}, l^*_{+\infty})$ an oriented universal cover.

A consequence of Proposition \ref{surfaceconvergence} is that the
lifts to $\mathcal{E}$ of the $\pi_i^{-1}(l^0_i\cup l^1_i \cup
l^{\infty}_i)$ converge for $d_H$ to the lift of
$\pi_{+\infty}^{-1}(\Lambda^+)$. We know that
$U=\pi_{+\infty}^{-1}(S_{+\infty}\smallsetminus \Lambda^+)\subset
\mathbb{H}^2$ is a disjoint union of (open) ideal triangles, of full
Lebesgue measure in $\mathbb{H}^2$. For any connected compact set
$K\subset U$, we have $K\cap \pi_i^{-1}(l^0_i\cup l^1_i \cup
l^{\infty}_i)=\emptyset$ for $i$ large enough, so $\varphi_{+\infty}$
is totally geodesic on $K$. Therefore, $\varphi_{+\infty}$ is totally
geodesic on each component of $U$. Since $\varphi_{+\infty}$ is
clearly $1$-Lipschitz, we can approximate any segment in
$\mathbb{H}^2\smallsetminus U$ by segments in $U$ to show that
$\varphi_{+\infty}$ is totally geodesic on each leaf of
$\mathbb{H}^2\smallsetminus U$. By Lemma 5.2.8 in \cite{lms},
$\varphi_{+\infty}$ sends rectifiable segments to rectifiable
segments of the same length, and is a pleated map.

\subsection{The map $\varphi_{+\infty}$ is a topological immersion} \llabel{patches}
Define $\mathcal{P}:=\pi_{+\infty}^{-1}(\Lambda^+)$, which contains
the pleating locus of $\varphi_{+\infty}$. To prove
$\varphi_{+\infty}$ is an immersion, it is enough to find a short
geodesic segment $m$ of $\mathbb{H}^2$, through the base point $\omega$, transverse to $\mathcal{P}$, and prove that $\varphi_{+\infty}$ is an immersion
on the union $\Upsilon$ of all strata (lines and complementary ideal triangles) of $\mathcal{P}$ crossed by $m$ (indeed,
$\pi_{+\infty}(\Upsilon)=S_{+\infty}$). Clearly, $\varphi_{+\infty}$ is already an immersion near any point of $\mathbb{H}^2\smallsetminus \mathcal{P}$. At other points, the key fact will be an ``equidistribution'' property of the $3$ pleating lines of the surface $S_i$, as $i$ goes to $+\infty$.

Choose a small $\mu_1>0$, and pick $k\in \mathbb{Z}$ large enough so
that $\phi_k^+\leq\mu_1$. Let $q$ be the rational opposite the Farey
edge $e_k$, on the same side as $\beta^-/\alpha^-$, so that
$\lambda^+(q)=\phi^+_k$. Let $m^{\mathcal{C}}_i$ be the simple
closed geodesic of slope $q$ in $S_i$ (equipped with the intrinsic
metric): for some superscript $*$ independent of $i$, the
Weierstrass point $\omega_i=\omega^*_i$ belongs to
$m^{\mathcal{C}}_i$. By Proposition \ref{surfaceconvergence} and
Hausdorff convergence, there exists $\mu_2>0$ such that for all
$i>k$, the angle between $m^{\mathcal{C}}_i$ and the pleating line
$l_i$ of $S_i$ at $\omega_i$ is at least $\mu_2$, and there exists
$\mu_3>0$ such that the segment $m_i$ of $m^{\mathcal{C}}_i$ of
length $2\mu_3$, centered at $\omega_i$, is embedded in $S_i$. The
ends of $m_i$ are at distance at least $\frac12 \mu_2\mu_3$ from the
pleating line $l_i$ in $S_i$. Observe that the simple closed curve
$m^{\mathcal{C}}_i$ meets the pleating edges of $S_i$ in a perfectly
equidistributed (Sturmian) order (as would a straight line on a Euclidean grid): therefore, the algebraic sum of the pleating angles crossed by any given subsegment of $m_i$ does not exceed $2\mu_1$. Finally, let $\kappa_i$ be a subsegment of the pleating line $l_i$, centered at $\omega_i$, of length $2$: if
$\mu_3$ is small enough, any pleating line $L$ met by $m_i$ makes an
angle at least $\frac12 \mu_2$ with $m_i$, and comes within $3\mu_3$
of both ends of $\kappa_i$ for the intrinsic metric of $S_i$.

\medskip
Arrange the developing map $\Phi:\widetilde{V}\rightarrow
\mathbb{H}^3$ in the upper half-space model so that the
$\varphi_i(\omega)$ lie on the line $0\infty$ at heights less than
$1$, and $\varphi_{+\infty}(l)$ is the oriented line from $-1$ to
$1$. Consider lifts $\widetilde{m}_i$ of the arcs $m_i$ through the
$\varphi_i(\omega)$. By the above (considering lifts of the
$\kappa_i$), if the $\mu$'s are small enough, any pleating line of
$\varphi_i(\mathbb{H}^2)$ met by $\widetilde{m}_i$ has its endpoints
within distance $1/2$ from $1$ and $-1$ in $\mathbb{C}$ (recall
$\varphi_i$ is $1$-Lipschitz). Following Subsection
\ref{definepleating}, let $c^{+1}_i$ (resp. $c^{-1}_i$) denote the
piecewise smooth curve defined by the transverse segment
$\widetilde{m}_i$ of $\varphi_i(\mathbb{H}^2)$ near $1$ (resp.
$-1$).

Let $\tau$ be a subsegment of $\widetilde{m}_i$ across an ideal
triangle of $\varphi_i(\mathbb{H}^2)$. Let $\tau'$ be the circle arc
contributed by $\tau$ to $c^{\pm 1}_i$. By the above, for some
universal $K_1>0$,
$$ \frac{\mu_2}{K_1} \leq \frac{\text{euclidean length
of }\tau'}{\text{hyperbolic length of }\tau} \leq K_1.$$ In
particular, the $c^{\pm 1}_i$ have length at most $2K_1\mu_3$. But
the regular curvature radii of $c^{\pm 1}_i$ are at least $\frac12$
(the corresponding circles come near $1$ and $-1$), so the total
regular curvature of $c^{\pm 1}_i$ is at most $4K_1\mu_3$. By the
above, the total singular curvature on any subinterval of $c^{\pm
1}_i$ is at most $2\mu_1$. If the $\mu$'s are small enough, it
follows that all tangent vectors of $c^{\pm 1}_i$ have complex
arguments within $[\pi/4, 3\pi/4]$. As a consequence, if $\tau$
ranges over the ideal triangles of $\varphi_i(\mathbb{H}^2)$ crossed
by $\widetilde{m}_i$, and $p:\mathbb{H}^3\rightarrow \mathbb{C}$ is
the vertical projection, the different $p(\tau)$ intersect only
along their edges: so $\varphi_i(\bigcup_{\tau}\tau)$ is an embedded
surface (which we can see, say, as the graph of a function from an
open set of $\mathbb{C}$ to $\mathbb{R}^+$). Moreover, define the
\emph{breadth} of $p(\tau)$ as the length of the segment
$p(\tau)\cap \sqrt{-1}\mathbb{R}$ (of the imaginary axis). Then, for some
universal $K_2>0$,
$$\frac{\text{breadth of }p(\tau)}{\text{hyperbolic length of }
\tau\cap \widetilde{m}_i} \geq \frac{\mu_2}{K_2}.$$

To conclude concerning $\varphi_{+\infty}$, define a lift
$m^{\mathcal{C}}$ through $\omega \in \mathbb{H}^2$ of the simple
closed geodesic of slope $q$ in $S_{+\infty}$, and a subsegment $m$
of $m^{\mathcal{C}}$, of length $2\mu_3$, centered at $\omega$. The
angle between $m$ and any pleating line of $\Lambda^+$ it encounters
is at least $\frac12 \mu_2$. Let $I$ be the (infinite) collection of
ideal triangles of $\pi_{+\infty}^{-1}(S_{+\infty}\smallsetminus
\Lambda^+)$ crossed by $m$. For $\tau\in I$, denote by $|\tau|$ the
length of $m\cap \tau$. By convergence in the Hausdorff metric,
$\varphi_{+\infty}(\tau)$ is approached by triangles of the
$\varphi_i(\mathbb{H}^2)$, of breadth at least
$\frac{|\tau|\mu_2}{2K_2}$: so $p(\varphi_{+\infty}(\tau))$ has nonzero
breadth. Injectivity follows: if $x,x' \in \overline{\bigcup_{\tau\in
I}\tau}$ do not belong to the same stratum of the lamination $\Lambda^+$, find $\tau$ separating
$x$ from $x'$ to prove that $\varphi_{+\infty}(x)\neq
\varphi_{+\infty}(x')$. By vertical projection to $\mathbb{C}$, we see
that $\varphi_{+\infty}(\mathbb{H}^2)$ is topologically immersed in
$\mathbb{H}^3$.

\subsection{Pleating measure of $\varphi_{+\infty}$} \llabel{pleatingphiinf}These arguments can be extended to prove that the pleating measure of
$\varphi_{+\infty}$, as defined in \ref{definepleating}, is the limit
of the pleating measures of the $\varphi_i$: the rectifiable
curves $c^{\pm 1}_{+\infty}$ defined near $\pm 1$ by
$\varphi_{+\infty}$ have lengths $\ell^{\pm 1}$, and for any
$\varepsilon >0$, there exists a finite disjoint union of circle
arcs $\gamma_s$ in $c^{\pm 1}_{+\infty}$ whose lengths add up to at
least $\ell^{\pm 1}-\varepsilon$ (moreover the direction of $c^{\pm 1}_{+\infty}$, like that of $c^{\pm 1}_i$, is everywhere within $\pi/4$ of the vertical axis). The $\gamma_s$ can be approached by (unions of) arcs of the $c^{\pm 1}_i$, and the regular curvature \emph{not} contributed by the $\gamma_s$ is bounded by $3\varepsilon$. It follows that the pleating of $\varphi_{+\infty}$ is the limit of the pleatings of the $\varphi_i$ (on any transverse arc, and therefore, on any simple closed curve): that pleating is simply $\lambda^+$. In particular, $\varphi_{+\infty}=\lim \varphi_{\nu(i)}$ is independent of the original subsequence $\nu$, and the $\varphi_i$ converge to a pleated map whose pleating is given by $\lambda^{\pm}$, as $i$ goes to $\pm \infty$.

\subsection{Completeness} The construction of \ref{patches} further
allows us to embed the universal cover $\widetilde{V}$ of
$V=\bigcup_{i\in \mathbb{Z}}\Delta_i$ into a (topological) manifold
with boundary $\widetilde{V}_{\partial}$, as follows. For each $x\in
\mathbb{H}^2$, consider a neighborhood $U_x$ of $x$ such that
$\varphi_{+\infty}$ is an embedding on $U_x$. Then
$\varphi_{+\infty}(U_x)$, which has a well-defined transverse
(``outward'') orientation, splits a small ball $B_x$ centered at
$\varphi_{+\infty}(x)$ into two (topological) hemispheres, which we
can call ``inner'' and ``outer'', referring to the transverse
orientation. The inner hemispheres $H_x$, for $x$ ranging over
$\mathbb{H}^2$, can be patched together to obtain a manifold with
boundary $H$. Without loss of generality, the balls $B_x$ can be
chosen small enough so that, by the construction of Section
\ref{patches}, each $H_x\smallsetminus \varphi_{+\infty}(U_x)$ is
identified with a subset of $\widetilde{V}$, embedded in
$\mathbb{H}^3$. Then, $H$ can further be patched to $\widetilde{V}$.
Since each $H_x$ is homeomorphic to
$\mathbb{R}^2\times\mathbb{R}^+$, the space
$\widetilde{V}_{\partial}=\widetilde{V}\cup H$ is a (possibly
non-complete) topological manifold with boundary.

\begin{proposition} \llabel{properlydiscontinuous}
The action of the fundamental group $\Gamma$ of the punctured torus $S$ on
$\widetilde{V}$ extends to a properly discontinuous action on
$\widetilde{V}_{\partial}$.
\end{proposition}
\begin{proof}
Consider the representation $\rho:\Gamma\rightarrow
\text{Isom}^+(\mathbb{H}^3)$ given by the developing map $\Phi$ from
(\ref{comdiagram}), and the representations
$\rho_n:\Gamma\rightarrow \text{Isom}^+(\mathbb{H}^2)$ which satisfy
$\varphi_n \circ \rho_n(g)=\rho(g) \circ \varphi_n$ for all $g\in
\Gamma$. By Lemma \ref{surfaceconvergence}, the $\rho_n$ converge to
some $\rho_{+\infty}$. Convergence of the $\varphi_n$ immediately
implies $\varphi_{+\infty} \circ \rho_{+\infty}(g)=\rho(g) \circ
\varphi_{+\infty}$. Therefore, the hemispheres $H_x$ can be chosen in
an equivariant fashion, and the action of $\Gamma$ on
$\widetilde{V}_{\partial}$ is well-defined. The action is already
properly discontinuous at every point $x$ of $\widetilde{V}$ (namely, as $g$ ranges over $\Gamma$, the $gx$ do not accumulate at $x$). But if
$\Gamma$ acts without fixed points and by isometries on a locally compact metric space $X$, the set of $x\in X$ such that the action is properly
discontinuous at $x$ is open (obviously) and closed: if $g_n x
\rightarrow x$ for some sequence $(g_n)$ of $\Gamma \smallsetminus
\{1\}$ and $U$ is a compact neighborhood of $x$, then $U$ contains a
ball of radius $\varepsilon$ centered at $x$, and whenever
$d(x,x')\leq \varepsilon/2$, the $g_n x'$ accumulate at some point
of $U$, so the action is not totally discontinuous at $x'$. The
Proposition follows by connectedness of $\widetilde{V}_{\partial}$.
\end{proof}
As a consequence, the space
$V_{\partial}:=\widetilde{V}_{\partial}/\Gamma$ is a (topological)
manifold with boundary, containing $V=\bigcup_{i\in
\mathbb{Z}}\Delta_i$; and $\partial V_{\partial}$ consists of two
pleated punctured tori (intrinsically isometric to $S_{+\infty}$ and
another surface $S_{-\infty}$), with pleatings $\lambda^+$ and
$\lambda^-$.

\begin{proposition}
The manifold with boundary $V_{\partial}$ is complete.
\end{proposition}
\begin{proof}
Consider the metric completion $\overline{V} \supset V_{\partial}$
and assume the inclusion is strict. Define a continuous function
$f:V_{\partial}\rightarrow \mathbb{R}^{>0}$ by
$f(x)=d(x,\overline{V} \smallsetminus V_{\partial})$. By assumption,
$\inf (f)=0$.

Consider the immersion of $\widetilde{V}_{\partial}$ into the upper
half-space model of $\mathbb{H}^3$ obtained by sending a lift of the
cusp to infinity (Figures \ref{frise} and
\ref{fourtriangles}-\ref{eighttriangles}). The image of $\partial
\widetilde{V}_{\partial}$ contains, in particular, vertical
half-planes (interrupted at some height above $\mathbb{C}$).
Therefore, any geodesic of $\widetilde{V}_{\partial}$ starting high
enough above $\mathbb{C}$ is defined for all times $t\leq 1$ (unless
it hits $\partial\widetilde{V}_{\partial}$). As a result, if
$H\subset V_{\partial}$ denotes a small enough open horoball
neighborhood of the cusp, we have $f\geq 1$ on $H$.

For $i\in \mathbb{Z}$, consider the compact set
$K_i:=S_i\smallsetminus H$ in $V_{\partial}$. There exists a ball
$B$ of $\mathbb{H}^2$ centered at the base point $\omega$, with radius independent
of $i$, such that $K_i\subset\pi_i(B)$: by convergence of $(\varphi_i)$,
the $\pi_i(B)$ converge metrically to a compact subset $K'$ of
$\partial V_{\partial}$, on which $f$ is positive. Therefore $f$ is bounded away from $0$ one some neighborhood $U$ of $K'$ in $V_{\partial}$, and $K_i\subset U$ for large enough $i$. So $f$ is bounded away from $0$ on $\bigcup_{i\in\mathbb{Z}} K_i$, and therefore on
$\bigcup_{i\in\mathbb{Z}} S_i$.

However, assume $\gamma(t)$ is a rectifiable $1$-Lipschitz arc of
$V_{\partial}$, defined for $t<M$, with no limit at $M$. For any
$\varepsilon>0$, the restriction $\gamma_{|[M-\varepsilon, M)}$
meets $V_{\partial}\smallsetminus
\partial V_{\partial}$ (because $\partial
V_{\partial}=S_{+\infty}\sqcup S_{-\infty}$ is complete), but then
$\gamma_{|[M-\varepsilon, M)}$ must meet $\Delta_i$ for an unbounded
set of indices $i$ (any finite union of tetrahedra is complete).
Therefore, we can find a sequence $t_n\rightarrow M$ such that
$\gamma(t_n)\in \bigcup_{i\in \mathbb{Z}} S_i$. Clearly,
$f(\gamma(t_n))\leq M-t_n$ which goes to $0$: a contradiction. So
$V_{\partial}$ is complete.
\end{proof}

\subsection{A quasifuchsian punctured-torus group}
The end of the argument is now quite standard: recall the complete
manifold with locally convex boundary $\widetilde{V}_{\partial}$,
which is a universal cover of $\overline{V}=V_{\partial}$. Given
distinct points $x,x'\in \widetilde{V}_{\partial}$, consider a
shortest possible path $\gamma$ from $x$ to $x'$. If $\gamma$ has an
interior point in $\partial\widetilde{V}_{\partial}$, by local
convexity, we must have $\gamma\subset \partial
\widetilde{V}_{\partial}$, and $\gamma$ is a geodesic segment of
$\partial\widetilde{V}_{\partial}$. If not, $\gamma$ is (the closure
of) a geodesic segment of $\widetilde{V}$. At any rate, the extended
developing map $\Phi:\widetilde{V}_{\partial}\rightarrow
\mathbb{H}^3$ is an embedding (it sends $\gamma$ to a segment with
distinct endpoints) and has a closed convex image $C$, endowed with
a properly discontinuous action of the fundamental group $\Gamma$ of
$S$ (Proposition \ref{properlydiscontinuous}). The action extends
properly discontinuously to $\mathbb{H}^3$ (this can be seen by
projecting any point of $\mathbb{H}^3$ to $C$). The manifold
$\mathbb{H}^3/\Gamma$ contains $\overline{V}\simeq C/\Gamma$, which has
the desired boundary pleatings $\lambda^{\pm}$. Clearly, $C$ is the
smallest closed, convex set containing all parabolic fixed points of
$\Gamma$; therefore $\overline{V}$ is the convex core, and $\Gamma$
is quasifuchsian, with pleating data $\lambda^{\pm}$. Theorem \ref{pureexistence} is proved.

\section{Application: the EPH theorem
}\llabel{sectioneph}

Let us quickly recall the correspondence between the horoballs of
${\mathbb H}^3$ and the vectors in the positive light cone of
Minkowski space. Endow ${\mathbb R}^4$ with the Lorentzian product
$\left \langle (x,y,z,t) | (x',y',z',t') \right \rangle :=
xx'+yy'+zz'-tt'$. Define $$X:=\{v=(x,y,z,t)\in{\mathbb R}^4 ~|~ t>0
\text{ and } \langle v | v \rangle = -1 \}.$$ Then
$\langle.|.\rangle$ restricts to a Riemannian metric on $X$ and
there is an isometry $X \simeq {\mathbb H}^3$, with
$\text{Isom}^+(X)$ a component of $SO_{3,1}({\mathbb R})$. We will
identify the point $(x,y,z,t)$ of $X$ with the point at Euclidean
height $\frac{1}{t+z}$ above the complex number $\frac{x+iy}{t+z}$
in the Poincar\'e upper half-space model. Under this convention, the
closed horoball $H_{d,\zeta}$ of Euclidean diameter $d$ centered at
$\zeta=\xi+i\eta\in {\mathbb C}$ in the half-space model corresponds
to  $\{v\in X ~|~ \langle v | v_{d,\zeta}\rangle\geq -1\}$, where
$v_{d,\zeta}= \frac{1}{d}(2\xi,2\eta,1-|\zeta|^2,1+|\zeta|^2)$. We
therefore identify $H_{d,\zeta}$ with the point $v_{d,\zeta}$ of the
isotropic cone. Similarly, the closed horoball $H_{h,\infty}$ of
points at Euclidean height no less than $h$ in the half-space model
corresponds to $\{v\in X ~|~ \langle v | v_{h,\infty} \rangle \geq
-1\}$ where $v_{h,\infty}=(0,0,-h,h)$, so we identify $H_{h,\infty}$
with $v_{h,\infty}$.

Consider the following objects: a complete oriented hyperbolic $3$-manifold $M$ with one cusp, a horoball neighborhood $H$ of the cusp, the universal
covering $\pi:{\mathbb H}^3\rightarrow M$, and the group $\Gamma
\subset \text{Isom}^+({\mathbb H}^3)$ of deck transformations of
$\pi$. Then $H$ lifts to a family of horoballs $(H_i)_{i\in I}$ in
${\mathbb H}^3$, corresponding to a family of isotropic vectors
$(v_i)_{i\in I}$ in Minkowski space. The closed convex hull $C$ of
$\{v_i\}_{i\in I}$ is $\Gamma$-invariant, and its boundary $\partial
C$ comes with a natural decomposition
$\widetilde{K}$ into polyhedral \emph{facets}. Let $U$ be the convex core of $M$, minus the
pleating locus. In \cite{comparing}, Akiyoshi and Sakuma extended
the Epstein-Penner convex hull construction to prove that
$\widetilde{K}$ defines an $H$-independent decomposition $K$ of $U$ into ideal hyperbolic polyhedra, typically tetrahedra, allowing for a few clearly defined types of degeneracies. In \cite{aswy1}, with Wada and Yamashita, they also conjectured

\begin{theorem} \llabel{eph}
If $M$ is quasifuchsian, homeorphic to the product of the punctured
torus with the real line, and has irrational pleating laminations of
slopes $\beta^+/\alpha^+$ and $\beta^-/\alpha^-$, then the
restriction of $K$ to the interior of the convex core of $M$ is
combinatorially the triangulation $(\Delta_i)_{i\in \mathbb{Z}}$
defined in Section \ref{sectionstrategy}.
\end{theorem}

\begin{proof}
It is known (see \cite{series} or Corollary \ref{plth} below) that
quasifuchsian punctured-torus groups are fully determined by their measured pleating laminations. It follows that the manifold constructed in Sections
\ref{sectionstrategy}-\ref{sectionextrinsicconvergence}, with the
same pleating data as $M$, is isometric to $M$. Recall the
triangulated space $\widetilde{V}\subset{\mathbb H}^3$ of
(\ref{comdiagram}), which is a universal cover of the union
$\bigcup_{i\in{\mathbb Z}}\Delta_i$ of all tetrahedra $\Delta_i$.
Given a horoball neighborhood $H$ of the cusp of $M$, consider its
lifts $(H_i)_{i\in I}$ in ${\mathbb H}^3$ and the corresponding
isotropic vectors $(v_i)_{i\in I}$ in Minkowski space ${\mathbb
R}^4$. Then each tetrahedron $\Delta$ of $\widetilde{V}$ has its
vertices at the centers of four horoballs $H_1,\dots,H_4$: we
introduce the convex hull $\tau_{\Delta}$ of $v_1,\dots,v_4$, and
define $D:=\bigcup_{\Delta} \tau_{\Delta}$. The central projection
to the hyperboloid $X$ with respect to the origin sends
$\tau_{\Delta}$ homeomorphically to the lift of $\Delta$ in $X$, so
the interiors of the $\tau_{\Delta}$ are pairwise disjoint, and each
$\tau_{\Delta}$ comes with a transverse orientation
$\overrightarrow{u}$ (given by any ray through $\Delta$ issued from
the origin). The theorem claims exactly that $D\subset\partial C$
(the inclusion is expected to be strict, since $\partial C$ also
contains faces projecting, say, to the boundary of the convex core
--- they are analyzed in detail in \cite{comparing}).

\medskip

Suppose that $D$ is locally convex, with $\overrightarrow{u}$
pointing inward. Define $I:=[1,+\infty)$ and, for any subset $Y$ of
a vector space, $IY:=\{ty~|t\in I, y\in Y\}$. Since the interior of
the convex core is convex, $ID$ is a convex set. Moreover, $ID$
contains all the $v_i$, so its closure $\overline{ID}$ contains
their convex hull $C$. Conversely, $D$ is clearly included in $C$
and it is easy to check that $IC\subset C$ (because $I\{v_i\}\subset
C$ for all $i$, see \cite{comparing}). So $ID \subset C$ and, by
closedness, $\overline{ID}=C$. Clearly, $D\subset \partial
\overline{ID}=\partial C$, as wished. So we only need to prove

\begin{lemma} \llabel{minkonwex}
The codimension-1 simplicial complex $D\subset {\mathbb R}^4$ is
locally convex ($\overrightarrow{u}$ pointing inward).
\end{lemma}

\emph{Proof ---} Consider adjacent ideal tetrahedra $\Delta,
\Delta'$ in ${\mathbb H}^3$ which are lifts from tetrahedra
$\Delta_{i-1}, \Delta_i$ of the manifold. We must prove that the
dihedral angle in $\mathbb{R}^4$ between $\tau_{\Delta}$ and $\tau_{\Delta'}$ points ``downward''. We will assume that the letter between $i-1$ and $i$
is an $L$ belonging to a subword $RL^nR$ of $\Omega$. In the link of the cusp, the pleated surface $S_i$ between $\Delta_{i-1}$ and $\Delta_i$ contributes
a broken line $(-1,\zeta,\zeta',1)$ in ${\mathbb C}$ together
with its iterated images under $u\mapsto u\pm 2$, as in Figure
\ref{minkowski} (we assume that the vertices $-1,1$ both belong to the
base segments of the Euclidean triangles just below and just above
the broken line, in the sense of Figure \ref{fourtriangles}). We use
the notation
\begin{eqnarray*}\zeta+1&=&\overrightarrow{a}=a\,e^{iA} \\ \zeta'-\zeta &
=&\overrightarrow{b}=b\,e^{iB} \\1-\zeta'&=&\overrightarrow{c}=c\,e^{iC}
\end{eqnarray*} (so far $A, B, C$
are only defined modulo $2\pi$). Above this broken line lives a lift of $S_i$ which admits as a deck transformation
$$f: u\mapsto 1+\frac{(\zeta+1)(\zeta'-1)}{u+1}~,$$ because $f(-1)=\infty$ ; $f(\infty)=1$ ; $f(\zeta)=\zeta'$. Therefore, $f(H_{1,\infty})=H_{|\zeta+1||\zeta'-1|,1}=H_{ac,1}$. Similarly, the following horoballs all belong to the same orbit: $$H_{1,\infty}~;~
H_{ac,-1}~;~ H_{ab,z}~;~ H_{bc,z'}~;~ H_{ac,1}~.$$ If $\zeta=\xi+\eta \sqrt{-1}$ and $\zeta'=\xi'+\eta'\sqrt{-1}$, the corresponding isotropic vectors in Minkowski space are respectively
\begin{equation} \begin{array}{lcccccccl}
v_{\infty}&=&&(&0,&0,&-1,&1&) \\ v_{-1}&=&\frac{1}{ac}&(&-2,&0,&0,&2&)\\
v_{\zeta}&=&\frac{1}{ab}&(&2\xi,&2\eta,&1-|\zeta|^2,&1+|\zeta|^2&)\\
v_{\zeta'}&=&\frac{1}{bc}&(&2\xi',&2\eta',&1-|\zeta'|^2,&1+|\zeta'|^2&) \\
v_{1}&=&\frac{1}{ac}&(&2,&0,&0,&2&). \end{array}\end{equation}

\begin{figure} [h!]\centering 
\psfrag{m1}{$(-1)$}
\psfrag{p1}{$(1)$}
\psfrag{z}{$(\zeta)$}
\psfrag{Z}{$(\zeta')$}
\psfrag{a}{$\overrightarrow{a}$}
\psfrag{b}{$\overrightarrow{b}$}
\psfrag{c}{$\overrightarrow{c}$}
\psfrag{ds}{$\overrightarrow{d_s}$}
\psfrag{dS}{$\overrightarrow{d_{s'}}$}
\psfrag{xi}{$x_i$}
\psfrag{yi}{$y_i$}
\psfrag{zi}{$\pi-w_i$}
\psfrag{xim}{$x_{i-1}$}
\psfrag{yim}{$y_{i-1}$}
\psfrag{zim}{$\pi-w_{i-1}$}
\ledessin{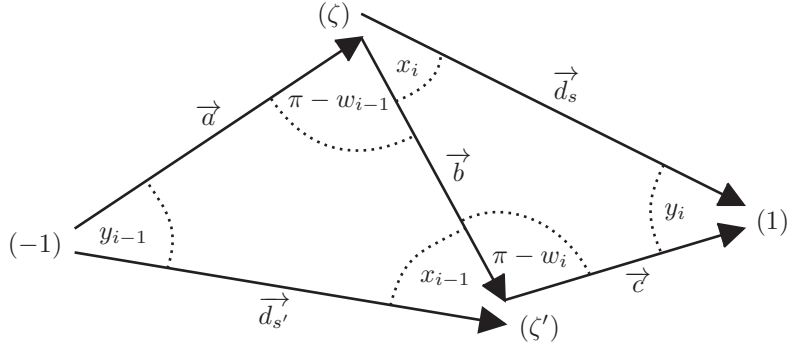}
\caption{Adjacent tetrahedra $\Delta_{i-1}, \Delta_i$ (cusp view).\llabel{minkowski}} \end{figure}

To prove that the dihedral angle at the codimension-2 face projecting to $(
\zeta \zeta' \infty)$ is convex, it is enough to show that if
$Pv_{\zeta}+Qv_{\zeta'}+Rv_{\infty}=\lambda v_1+(1-\lambda)v_{-1}$ then
$P+Q+R>1$ (moreover, this will in fact take care of \emph{all}
codimension-2 faces of the simplicial complex $D$). One easily finds the unique solution
$$P=\frac{-b\eta'}{c(\eta-\eta')}~;~Q=\frac{b\eta}{a(\eta-\eta')}~;~R= \frac{\eta(1-|\zeta'|^2)-\eta'(1-|\zeta|^2)}{ac(\eta-\eta')}$$
hence $$P+Q+R=1+\frac{Z}{ac(\eta-\eta')}~~\text{ where }
Z=bc\eta-ab\eta'+\eta(1-|\zeta'|^2)-\eta'(1-|\zeta|^2)+ac(\eta'-\eta).$$ Observe that $\eta>\eta'$
because the triangles $-1\zeta'\zeta$ and $1\zeta\zeta'$ are counterclockwise
oriented. So it is enough to prove that $Z>0$. Endow ${\mathbb
C}\simeq {\mathbb R}^2$ with the usual scalar product $\bullet$ and
observe that $1-|\zeta|^2=\overrightarrow{a}\bullet
(\overrightarrow{b}+\overrightarrow{c})$ and
$1-|\zeta'|^2=(\overrightarrow{a}+\overrightarrow{b})\bullet\overrightarrow{c}$.
Hence \begin{eqnarray*} Z&=&
\eta(bc+\overrightarrow{b}\bullet\overrightarrow{c})-\eta'(ab+\overrightarrow{a}\bullet\overrightarrow{b})+
(\eta'-\eta)(ac-\overrightarrow{a}\bullet \overrightarrow{c}) \\ 
&=& abc\left [\frac{\eta}{a}(1+\cos(B-C))-\frac{\eta'}{c}(1+\cos(A-B))+\frac{\eta'-\eta}{b}(1-\cos(A-C))\right ] \\ &=& abc[\sin A(1+\cos(B\!-\!C)) + \sin C(1+\cos(A\!-\!B)) + \sin B(1-\cos(A\!-\!C))] \\ &=&
4abc\,\sin\frac{A+C}{2}\cos\frac{A-B}{2}\cos\frac{B-C}{2}.
\end{eqnarray*} by standard trigonometric formulae. Observe that the
last expression is a well-defined function of $A,B,C \in \mathbb{R}/2\pi\mathbb{Z}$ (although each factor is defined only up to sign). Next, however, we shall
carefully pick representatives of $A,B,C$ in ${\mathbb R}$.

\begin{observation} \llabel{argumentnormal}
There exists a unique triple of representatives $(A,B,C) \in \mathbb{R}^3$ such that $\{A,B,C\}\subset J$ where $J$ is an open interval of length less than $\pi$, containing $0$. (It is easy to see that the broken line $(\dots,-1,\zeta,\zeta',1,\dots)$ has no self-intersection if and only if there exists an open half-plane $H$ containing the vectors $\overrightarrow{a}, \overrightarrow{b}, \overrightarrow{c}$: the observation follows since $\overrightarrow{a}+\overrightarrow{b}+\overrightarrow{c}=2$ belongs to $H$). \end{observation}
We pick the representatives $A,B,C \in \mathbb{R}$ given by Observation \ref{argumentnormal}, and do the same for the broken line contributed by each surface $S_i$ for $i\in\mathbb{Z}$ (each broken line is oriented from $-\infty$ to $+\infty$, so each edge $s$ of the cusp link inherits an orientation, as in Figure \ref{minkowski}: we write its complex coordinate
$\overrightarrow{d_s}=d_s\, e^{D_s}$). By Observation \ref{argumentnormal}, all the complex arguments $D_s$ are in $(-\pi,\pi)$. But since $0\leq w \leq \pi$, the existence of the interval $J$ also implies (in the notation of Figure \ref{minkowski}):
$$A=B+w_{i-1}~~~\text{ and }~~~C=B+w_i$$
hence $B=\inf J$ and $B\in (-\pi,0)$. It follows that $B+x_i \in (-\pi,\pi)$, hence $$D_s=B+x_i~~~\text{ and similarly }~~~ D_{s'}=B+x_{i-1}.$$
In other words (by transitivity), for any two edges $s_1,s_2$ of the cusp link, the difference of arguments $D_{s_1}-D_{s_2}\in\mathbb{R}$ can be read off ``naively'' as a linear combination of the $\{w_i\}_{i\in\mathbb{Z}}$, with no multiple of $2\pi$ added.


Now that all arguments are fixed in the real interval $(-\pi, \pi)$,
we can see that $\cos\frac{A-B}{2}=\cos (w_{i-1}/2)$ and
$\cos\frac{B-C}{2}=\cos (w_i/2)$ are positive. So to prove $Z>0$ it
remains to see that $0<A+C<2\pi$. This in turn follows from a small
variant of Lemma 16 in \cite{mapomme}, which we prove now (it implies an empirical observation of J{\o}rgensen which was Conjecture 8.6 in \cite{aswy1}).
\begin{proposition} \llabel{tiltlemma}
With the above notation, one has $0<A+C<2\pi$.
\end{proposition}
\begin{figure}[h!] \centering 
\psfrag{a}{$\overrightarrow{a}$}
\psfrag{c}{$\overrightarrow{c}$}
\psfrag{ceq}{$\left (\!\! \begin{array}{l} \overrightarrow{a}=a\,e^{iA} \\ \overrightarrow{c}=c\,e^{iC} \end{array} \!\! \right )$}
\psfrag{m1}{$-1$}
\psfrag{p1}{$1$}
\psfrag{zm1}{$\zeta_{-1}$}
\psfrag{z0}{$\zeta_0$}
\psfrag{z1}{$\zeta_1$}
\psfrag{znm}{$\zeta_{n-1}$}
\psfrag{zn}{$\zeta_n$}
\psfrag{znp}{$\zeta_{n+1}$}
\psfrag{ze0}{$\xi_0$}
\psfrag{zen}{$\xi_n$}
\psfrag{zep}{$\xi_{n+1}$}
\ledessin{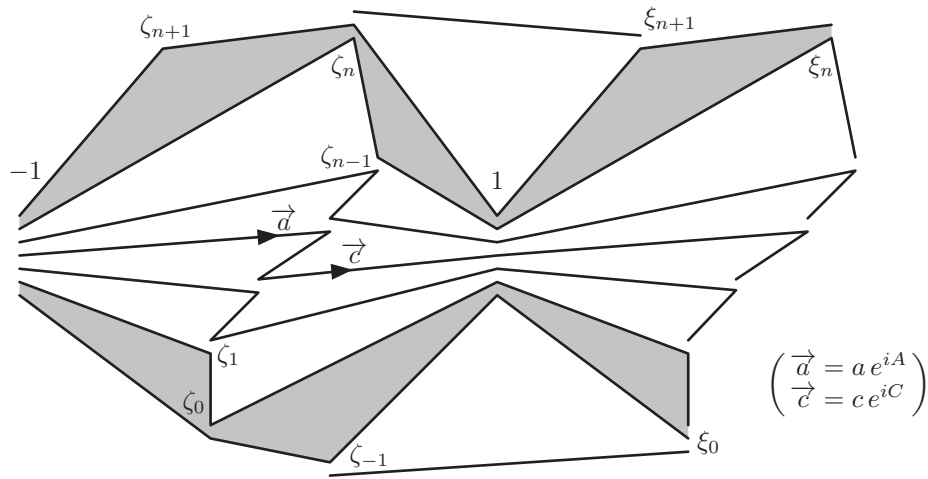}
\caption{\llabel{manytriangles} Cusp view: a subword $RL^nR$, bounded by hinge indices $0$ and $n$.} \end{figure}

Consider the maximal subword $RL^nR$, with the broken line
$(-1,\zeta,\zeta',1)$ corresponding to some $L$. There is in fact a sequence
of points $(\zeta_0, \zeta_1, \dots,\zeta_n)$ in ${\mathbb C}$ such that the
broken line corresponding to the $j$-th letter $L$ is
$(-1,\zeta_j,\zeta_{j-1},1)$ for all $1\leq j \leq n$. The broken line
corresponding to the initial (resp. final) $R$ is $(-1,\zeta_0,\zeta_{-1},1)$
for some $\zeta_{-1}$ (resp. $(-1,\zeta_{n+1},\zeta_n,1)$ for some $\zeta_{n+1}$). There exists $1\leq k \leq n$ such that $(\zeta,\zeta')=(\zeta_k,\zeta_{k-1})$.

Set $\xi_j:=\zeta_j+2$ for all $j$. By construction, the rays issued from $1$ through $\zeta_{-1},\zeta_0,\dots,\zeta_{n-1},\zeta_n,\xi_{n+1},\xi_n,\dots,\xi_1,\xi_0,\zeta_{-1}$ (in that cyclic order) divide ${\mathbb C}$ clockwise into salient angular sectors of sum $2\pi$ (these rays realize the link of the
vertex $1$ in the cusp triangulation, see Figure \ref{manytriangles}).
Comparing the angles at $1$ and at $-1$, we see immediately that
\begin{equation} \llabel{convert} A+C=D_{-1 \zeta_k}+D_{\zeta_{k-1} 1}=D_{-1
\zeta_{k+1}}+D_{\zeta_k 1}=\dots=D_{-1 \zeta_j}+D_{\zeta_{j-1} 1} \end{equation}
(for all $0\leq j \leq n+1$). Since the triangles $\zeta_0 \zeta_{-1}1$
and $\xi_0 1\zeta_{-1}$ are counterclockwise oriented, we have
${\rm Im}(\zeta_{-1})<{\rm Im}(1)$ i.e. ${\rm Im}(\zeta_{-1})<0$.
Similarly, $\zeta_n 1 \xi_{n+1}$ and $\xi_n \xi_{n+1}1$ are
counterclockwise oriented, so ${\rm Im}(\xi_{n+1})>0$. So there exists
$0\leq j \leq n+1$ such that ${\rm Im}(\zeta_{j-1})< 0$ and ${\rm
Im}(\zeta_j)\geq 0$. This implies $D_{-1 \zeta_j} \in (0,\pi)$ and
$D_{\zeta_{j-1} 1}\in[0,\pi)$. By (\ref{convert}), this implies $A+C\in
(0,2\pi)$. Theorem \ref{eph} is proved. \end{proof}

In \cite{mapomme}, we studied punctured-torus bundles over the circle by indexing the tetrahedra $\Delta_i$ in $\mathbb{Z}/m\mathbb{Z}$ (instead of $\mathbb{Z}$): note that the proof of Theorem \ref{ephgeneral} applies without alteration to that context.

\section{Generalizations} \llabel{sectionextensions}
In this section, we extend all the previous results to punctured-torus groups with rational pleatings and/or infinite ends. The general theorem is as follows.
\begin{theorem} \llabel{ephhyper}
Let $\lambda^+\neq \lambda^-$ be nonzero, possibly projective (if irrational), measured laminations on the punctured torus $S$, and let $s(\lambda^{\pm})\in \mathbb{P}^1\mathbb{R}$ be the slope of $\lambda^{\pm}$. If $\lambda^{\pm}$ is a closed leaf, denote its weight by $|\lambda^{\pm}|$ and assume $|\lambda^{\pm}|\leq\pi$. There exists a punctured-torus group $\Gamma$ with ending and/or pleating laminations $\lambda^{\pm}$, and the open convex core $V$ of $\mathbb{H}^3/\Gamma$ has an ideal decomposition $\mathcal{D}$ into polyhedral cells (of positive volume) whose combinatorics are given by $\lambda^{\pm}$ in the following sense: if $\Lambda$ is the line from $s(\lambda^-)$ to $s(\lambda^+)$ across the Farey diagram in $\mathbb{H}^2$, then
\begin{enumerate}
\item If $s(\lambda^+)$ and $s(\lambda^-)$ are irrational, $\mathcal{D}$ consists of ideal tetrahedra $(\Delta_i)_{i\in\mathbb{Z}}$ in natural bijection with the Farey edges crossed by $\Lambda$, as in Section \ref{sectionstrategy}.
\item If only $s(\lambda^+)$ is rational and $|\lambda^+|<\pi$, then $\mathcal{D}$ has one ideal tetrahedron per Farey edge crossed by $D$, and one cell $T$ homeomorphic to a solid torus: $\partial_1=\overline{\partial T \cap \partial V}$ is a punctured torus pleated along a simple closed curve of slope $s(\lambda^+)$, and $\partial_2=\overline{\partial T \smallsetminus \partial V}$ is a punctured torus pleated along the ideal triangulation associated to the Farey triangle with vertex $s(\lambda^+)$ crossed by $\Lambda$. Finally, $\partial_1\cap\partial_2$ is a line from the puncture to itself of slope $s(\lambda^+)$. See the left panel of Figure \ref{toricpieces}.
\item If only $s(\lambda^+)$ is rational and $|\lambda^+|=\pi$, all the statements of the previous case apply, except that $\partial_1=\overline{\partial T \cap \partial V}$ becomes a thrice-punctured sphere (the simple closed curve of slope $s(\lambda^+)$ has been ``pinched'' to become a cusp): see the right panel of Figure \ref{toricpieces}.
\item If only $s(\lambda^-)$ is rational, the situation is similar to the two previous cases, exchanging $\lambda^-$ and $\lambda^+$.
\item If $s(\lambda^+),s(\lambda^-)$ are rationals but not Farey neighbors, the situation is again similar, with two solid torus cells $T^+$ and $T^-$.
\item If $s(\lambda^+), s(\lambda^-)$ are Farey neighbors, $\mathcal{D}$ only consists of two solid tori $T^+$ and $T^-$ as above, glued along a punctured torus $S$ pleated along only two lines, of slopes $s(\lambda^+)$ and $s(\lambda^-)$.
\end{enumerate}
Moreover, $\mathcal{D}$ agrees with the geometrically canonical decomposition $\mathcal{D}^G$ of the open convex core given by the Epstein-Penner convex hull construction.
\end{theorem}
Note that the combinatorics of $\mathcal{D}^G$ do not depend on the nature (finite or infinite) of the ends of $\mathbb{H}^3/\Gamma$. At this point, we have treated the case of two irrational pleatings (finite ends). We proceed to prove the remaining cases of Theorem \ref{ephhyper}.

\begin{figure}[h!] \centering 
\psfrag{tlp}{$\theta<\pi$}
\psfrag{tp}{$\theta=\pi$}
\ledessin{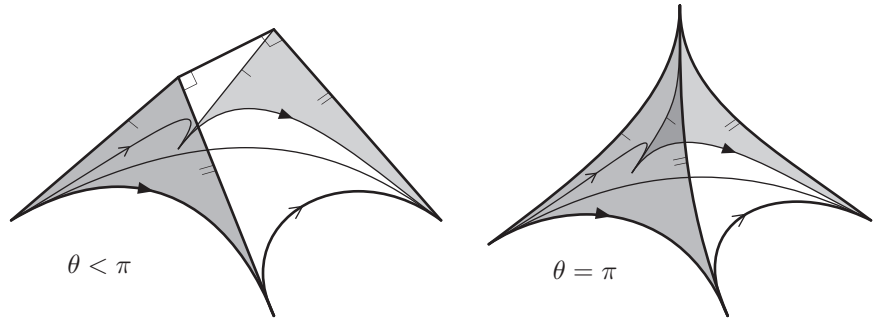}
\caption{Toric cells: the exterior dihedral angle $\theta$ is the weight $|\lambda^{\pm}|$. Shaded faces are identified. \llabel{toricpieces}} \end{figure}

\subsection{One rational pleating} \llabel{onerationalpleating}
We focus on the case of two finite ends, with only $\beta^+/\alpha^+$ rational.
We can choose to end the word $\Omega\in \{R,L\}^{\mathbb{Z}}$ with an infinite suffix $LRR...R...$ (or of $RLL...L...$: that is an arbitrary choice) and proceed from
Section \ref{sectionstrategy} onward. We shall assume that $i=0$ is
the greatest hinge index. Sections
\ref{sectionangles} through \ref{subsectionnaturalconstraint} are unchanged:
the sequence $(w_i)$ is just concave (thus convergent and non-decreasing) on
$\mathbb{N}$. In Subsection \ref{subsectionstudyphi}, we find that
the sequence $(\phi^+_i)$ is constant on $\mathbb{N}$, equal to some
$\theta>0$. By (\ref{phiphi}), the number $\theta$ is the weight of the rational lamination $\lambda^+$, so we assume $\theta<\pi$.
Section \ref{sectionhyperbolicvolume} goes through
essentially unchanged: by the computation of Sublemma \ref{boundhops} (and with the same notation), the sum of the volumes of all tetrahedra $\Delta_i$ for $i\geq 2^n$ is at most
$$\sum_{k>n} \Sigma_{2^{k-1}}^{2^k}\leq \sum_{k>n}
2^{-k}[1+(2k-1)\log 2]=O(2^{-n/2}).$$ Therefore the volume
functional $\mathcal{V}$ is bounded, continuous for the product
topology, and concave.
We can find a maximizer $w$ of $\mathcal{V}$, and it still satisfies Propositions \ref{tricot} and \ref{hingefragile}: in particular, all tetrahedra $\Delta_i$ for $i>0$ have positive angles.

Something new is required in Section \ref{sectionbehaviorofwi}: we
must prove that $\lim_{i\rightarrow +\infty}w_i=\theta$.

The $\{w_i\}_{i\geq 2}$ contribute only to the angles of the $\{\Delta_i\}_{i\geq 1}$, which are positive: so the volume $\mathcal{V}$ is critical with respect to each $w_i$ for $i\geq 2$. By Sublemma 6 of \cite{mapomme}, it follows that the cusp triangles of the $\{\Delta_i\}_{i\geq 1}$ fit together correctly and can be drawn in the Euclidean plane $\mathbb{C}$. More precisely, there exists a sequence of complex numbers $(\zeta_i)_{i\geq 0}$ such that the triangles contributed by $\Delta_i$ have vertices at $(-1,\zeta_i,\zeta_{i-1})$ and $(1,\zeta_i,\zeta_{i+1})$ (Figure \ref{nautilus}, left). These triangles being similar, we have $(\zeta_i+1)(\zeta_{i+1}-1)=(\zeta_{i-1}+1)(\zeta_i-1)=...=(\zeta_0+1)(\zeta_1-1)$, hence
$$\zeta_{i+1}=\frac{\zeta_i+\kappa}{\zeta_i+1}=:\varphi(\zeta_i)$$ for some complex number $\kappa\neq 1$ independent of $i$. Observe that the complex length $\ell$ of the hyperbolic isometry (extending) $\varphi$ satisfies $\cosh \ell = \frac{1+\kappa}{1-\kappa}$. In the end, $\varphi$ will be a lift of the loop along the rational pleating line of the convex core.

\begin{proposition}
The number $\kappa$ lies in the real interval $(0,1)$.
\end{proposition}
\begin{proof} Denote by $Z_i$ the periodic broken line $(\dots,-1,\zeta_{i-1},\zeta_i,1,2+\zeta_{i-1},2+\zeta_i,\dots)$ contributed by the pleated surface $S_i$. First, $\kappa$ is real: if not, the $\zeta_i$ have a limit in $\mathbb{C} \smallsetminus \mathbb{R}$ (a square root of $\kappa$), so $w_i-w_{i-1}$ (the angle of $Z_i$ at $1$) cannot go to $0$. If $\kappa<0$, then $\varphi$ is a pure rotation: the $\zeta_i$ all belong to a circle of the lower half-plane, which also contradicts $w_i-w_{i-1}\rightarrow 0$. If $\kappa>1$, then $\varphi$ is a gliding axial symmetry: the $\zeta_i$ go to $\pm\sqrt{\kappa}$ and belong alternatively to the upper and lower half-plane, so $Z_i$ has self-intersection for large $i$. If $\kappa=0$, then $\varphi$ is a parabolic transformation fixing $0$: the $\zeta_i$ go to $0$ along a circle tangent to $\mathbb{R}$, and one can see that $w_i=\widehat{\zeta_i \zeta_{i+1}1}$ goes to $\pi>\lim_{+\infty}\phi^+=\theta$. The only remaining possibility is $\kappa \in (0,1)$, where $\varphi$ is a pure translation.
\end{proof}

Thus, the $\zeta_i$ lie on a circle arc $C$ of the lower half-plane which meets the real line at $\pm \sqrt{\kappa}$, and $\lim_{+\infty} \zeta_i=\sqrt{\kappa}$. We denote by $\theta^*$ the angle between $C$ and the segment $[-\sqrt{\kappa},\sqrt{\kappa}]$ (more precisely, the angle between their half-tangents at $\sqrt{\kappa}$). It is easy to see that $$\theta^*=\underset{i\rightarrow +\infty}{\lim} w_i$$ (see Figure \ref{nautilus}).
Hence, $\theta^*\leq \theta$. 

\begin{figure}[h!] \centering 
\psfrag{m1}{$-1$}
\psfrag{p1}{$1$}
\psfrag{msqk}{$-\sqrt{\kappa}$}
\psfrag{sqk}{$\sqrt{\kappa}$}
\psfrag{C}{$C$}
\psfrag{z0}{$\zeta_0$}
\psfrag{z1}{$\zeta_1$}
\psfrag{z2}{$\zeta_2$}
\psfrag{fz0}{$f(\zeta_0)$}
\psfrag{fz1}{$f(\zeta_1)$}
\psfrag{fz2}{$f(\zeta_2)$}
\psfrag{0}{$0$}
\psfrag{r}{$\rho$}
\ledessin{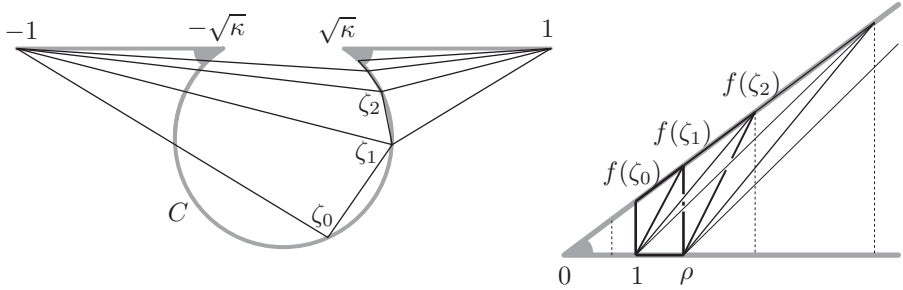}
\caption{All marked angles (grey) are $\pi-\theta^*$, and $\zeta_{j+1}=\varphi(\zeta_j)=\frac{\zeta_j+\kappa}{\zeta_j+1}$. \llabel{nautilus}} \end{figure}

\begin{proposition} One has $\theta^*=\theta$. \llabel{behaviorofwirational} \end{proposition}
\begin{proof}
First, it is easy to check that any data $0\leq w_0<w_1<\theta^*< \pi$ smoothly determines a \emph{unique} pair of complex numbers $\zeta_0,\zeta_1$ such that: 
\begin{enumerate}
\item the broken line $(-1,\zeta_0,\zeta_1,1,2+\zeta_0,2+\zeta_1,\dots)$ has angles $(w_0,-w_1,w_1-w_0)$ as in (\ref{pleatingangles}) above;
\item the number $\kappa$ such that $\zeta_1=\frac{\zeta_0+\kappa}{\zeta_0+1}$ lies in $(0,1)$;
\item the circle through $\zeta_0$ and $\zeta_1$ centered on the imaginary axis intersects the real axis at an angle $\theta^*$. \end{enumerate}

These $\zeta_0,\zeta_1$ in turn define all $\{\zeta_j\}_{j\geq 2}$ \emph{via} $\zeta_{j+1}=\varphi(\zeta_j)$, and we can read off the angle $w_j=\widehat{1 \zeta_{j+1} \zeta_j}\leq \theta^*$ and construct the associated ideal tetrahedron $\Delta_j$. In what follows, we investigate the shape of the space $U:=\bigcup_{j\geq 1} \Delta_j$, whose boundary (the punctured torus $S_1$) has pleating angles $(w_0,-w_1,w_1-w_0)$ (see (\ref{pleatingangles}) above).

Define $f(\zeta):=\frac{\zeta+\sqrt{\kappa}}{\zeta-\sqrt{\kappa}}$, so that $f(\varphi(\zeta))=\rho f(\zeta)$ where $\rho:=\frac{1+\sqrt{\kappa}}{1-\sqrt{\kappa}}$. The convex hull of $(\infty,1,\zeta_j,\zeta_{j+1})$ is isometric to the tetrahedron $\Delta_j$: pushing forward by $f$, we obtain a tetrahedron $\Delta'_j$, isometric to $\Delta_j$, with vertices $(1,\rho, \rho^j f(\zeta_0),\rho^{j+1} f(\zeta_0))$ (Figure \ref{nautilus}, right). Moreover, all the $f(\zeta_j)=\rho^j f(\zeta_0)$ lie on the half-line $e^{i(\pi-\theta^*)}\mathbb{R}^+$. By reasoning in a fundamental domain of the loxodromy $\Phi: u\mapsto \rho u$, it is then easy to see that $U$ has the same volume as $D:= D_1 \cup D_2$, where $D_1,D_2$ are ideal tetrahedra of vertices $(\infty, f(\zeta_1),1,\rho)$ and $(\infty,f(\zeta_1),1,f(\zeta_0))$ respectively. Moreover, $\Phi$ identifies the faces $(\infty,1,f(\zeta_0))$ and $(\infty,\rho,f(\zeta_1))$ of $D$, so that $D/\Phi$ is a manifold with polyhedral boundary, homeomorphic to a solid torus, with interior dihedral angles $$(\pi-w_1,\pi+w_0,\frac{w_1-w_0}{2},\frac{w_1-w_0}{2},\pi-\theta^*).$$ The edge of $D/\Phi$ with dihedral angle $\pi-\theta^*$ is a simple closed curve of length $\log \rho$, toward which $D_1,D_2$ spiral. A picture of $D/\Phi$ is obtained by replacing $\theta$ with $\theta^*$ in the left panel of Figure \ref{toricpieces}.

Using the smooth dependence on $\theta^*$, the Schl\"afli volume formula then gives $d\mathcal{V}(D/\Phi)/d\theta^*=\frac{1}{2}\log \rho>0$: so $D/\Phi$ (and therefore $U$) has largest volume when $\theta^*$ is largest. Regard the $\{w_j\}_{j\leq 1}$ as fixed, and the $\{w_j\}_{j\geq 2}$ as unknowns: then $(w_2, w_3,\dots)$ is clearly the solution to the maximization problem for the volume of $U$ (with fixed pleating angles on $S_1$). Therefore, the $w_i$ will choose the largest possible limit $\theta^*$ at $+\infty$, namely $\theta^*=\theta$.
\end{proof}

The above proof does more than determining $\lim_{+\infty} w_j$: as in Sections \ref{sectionbehaviorofwi}-\ref{sectionextrinsicconvergence}, it gives a full description of $U=\bigcup_{j\geq 1} \Delta_j$ and of its boundary (whose pleating turns out to be $\lambda^+$). Here is, however, an important observation:
\begin{observation} \llabel{depli}
Proposition \ref{tiltlemma} above, ``$A+C>0$'', does not hold for the family of pleated surfaces $(S_i)_{i\geq 0}$. Instead, we have $A+C=0$.
This simply means that $\overrightarrow{-1,\zeta_i}$ and $\overrightarrow{\zeta_{i+1},1}$ make opposite angles with the real line. Indeed, $(\zeta_i+1)(1- \zeta_{i+1})=1-\kappa^2$ is a positive real.\end{observation}

We can now establish 
\begin{proposition} All the (strict) inequalities of (\ref{doubleve}) are true.\llabel{winteriorrational} \end{proposition}
\begin{proof}
This is Proposition \ref{winterior} (in the new context where $\beta^+/\alpha^+$ is rational). The proof is the same, with the following caveat: in ruling out $w_j=0$ for $j$ hinge, we call upon \cite{mapomme} (especially Lemma 16 and the argument of Section 9 there). The strategy is to assume $w_j=0$, and then perturb $w$ to a well-chosen $w^{\varepsilon}$ so as to make the volume increase: $\partial\mathcal{V}/\partial\varepsilon>0$. The latter inequality holds essentially because the inequality of Proposition \ref{tiltlemma} (``$A+C>0$'') is true, both in the $R^n$-word preceding $j$ and in the $L^m$-word following $j$ (\footnote{In Section 9 of \cite{mapomme}, the inequality ``$A+C>0$'' is formulated in terms of lengths and takes the guise ``$Q<P+T$''.}). More exactly, $\partial\mathcal{V}/\partial\varepsilon$ will be positive when at least one of the two instances of ``$A+C\geq 0$'' is strict. But it is always strict, except in a single case (the infinite suffix $LRR...R...$): so we can conclude.
\end{proof}

As a result, all the tetrahedra $\Delta_i$ have positive angles and fit together correctly: Sections \ref{sectionbehaviorofwi}-\ref{sectionextrinsicconvergence} carry through for the $\lambda^-$-end, and $\bigcup_{i\in\mathbb{Z}}\Delta_i$ is the open convex core of a quasifuchsian punctured-torus group, with the prescribed pleatings $\lambda^{\pm}$.

The results of Section \ref{sectioneph} extend readily: the only modification is that tetrahedra $(\Delta_i)_{i\geq 0}$ lift to a family of \emph{coplanar} cells in Minkowski space, because the key inequality of Proposition \ref{tiltlemma} has become an equality. Therefore the geometrically canonical decomposition of the convex core contains the non-contractible cell $D/\Phi$ (Figure \ref{toricpieces}, left).

\subsection{Two rational pleatings} \llabel{tworationalpleatings}
When both $\beta^+/\alpha^+$ and $\beta^-/\alpha^-$ are rational, the encoding word $\Omega\in\{R,L\}^{\mathbb{Z}}$ can be chosen with an infinite prefix $...R...RRL$ and an infinite suffix $LRR...R...$ We apply the same argument as above to both ends simultaneously. Again, Proposition \ref{winteriorrational} holds because no hinge index $j$ belongs both to the prefix and to the suffix.

If the rational pleating slopes $s(\lambda^{\pm})$ are not Farey neighbors, we can convert the prefix to $...LLR$ and/or the suffix to $RLL...$, obtaining different triangulations of the convex core of the same quasifuchsian group.

If $s(\lambda^+),s(\lambda^-)$ are Farey neighbors, then $\Omega=...RRLRR...$ (observe that prefix and suffix overlap, so the two word conversions do \emph{not} commute). If the indices before and after the $L$ are $0$ and $1$, we obtain $w_0=w_1$ by applying Observation \ref{depli} to prefix and suffix. In other words, the points $a,b,c$ in Figure \ref{fareyneighbors} (left) are collinear. It is therefore possible to triangulate the same convex core according to a word $\Omega=...LLLRRR...$, provided that we allow the hinge tetrahedron $\Delta_0$ to become flat (black in Figure \ref{fareyneighbors}, right). In any case, the sequence $(w_i)_{i\in\mathbb{Z}}$ maximizes the total volume.

\begin{figure}[h!] \centering 
\psfrag{rlr}{$...RRLRR...$}
\psfrag{lr}{$...LLLRRR...$}
\psfrag{a}{$a$}
\psfrag{b}{$b$}
\psfrag{c}{$c$}
\ledessin{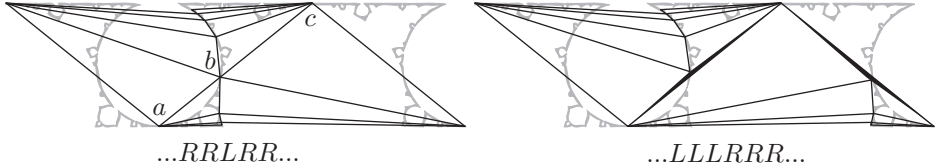}
\caption{Two triangulations seen against the same limit set. \llabel{fareyneighbors}} \end{figure}

\subsection{Pinching}
The case where one (or both) of the pleatings $\lambda^{\pm}$ has weight $\pi$ is a straightforward limit case of Subsections \ref{onerationalpleating} and \ref{tworationalpleatings} (the term ``pinching'' refers to the fact that the pleating curve becomes shorter and shorter, and eventually turns into a cusp as the pleating angle reaches $\pi$).

Suppose $|\lambda^+|=\pi$ (note that the conditions of (\ref{doubleve}) involving $\phi^+$ become vacuous, because $\phi^+\geq\pi$). Subsection \ref{onerationalpleating} carries through with $\theta=\pi$ and $\kappa=0$. The circle arc $C$ of Figure \ref{nautilus} becomes a full circle, tangent to $\mathbb{R}$ at $0$. The analysis of Proposition \ref{behaviorofwirational} (existence and uniqueness of $\zeta_0, \zeta_1\in\mathbb{C}$) extends smoothly to $\theta^*=\pi$. One then finds that the tetrahedron with vertices $(\infty,1,\zeta_j,\zeta_{j+1})$ is sent by $f:u\mapsto 1/u$ to a tetrahedron $\Delta'_j$ of $\mathbb{H}^3$, isometric to $\Delta_j$, with vertices $(0,1,\tau+j,\tau+j+1)$ for a certain $\tau \in \mathbb{C}$ independent of $j$ (in the right panel of Figure \ref{nautilus}, the grey bounding rays are replaced by parallel lines). It is straightforward to check that $\bigcup_{j\geq 1}\Delta_j$ is the solid torus pictured in Figure \ref{toricpieces} (right), and the Schl\"afli formula again implies that $\theta^*=\pi$ realizes the maximum volume.

After Proposition \ref{behaviorofwirational}, the argument is unchanged.

\subsection{Infinite ends}
By \cite{minsky}, quasifuchsian groups are dense in the set of all discrete, faithful, type-preserving representations of $\pi_1(S)\rightarrow \text{Isom}^+(\mathbb{H}^3)$. In fact, a geometrically infinite end of such a representation comes with an \emph{ending lamination}, namely an irrational projective measured lamination which should be thought of as an ``infinitely strong pleating''. In this section we will assume $\beta^+/\alpha^+ \notin \mathbb{P}^1 \mathbb{Q}$ and consider $w^T$, the solution of the volume maximization problem for $(\phi^-,T\phi^+)$, where $T>0$ (namely, $w$ is subject to conditions (\ref{doubleve}), where $\phi^+$ is replaced by $T \phi^+$; and Table (\ref{xiyi}) still expresses $x_i, y_i, z_i$).

In a discrete, type-preserving representation of a surface group, each measured lamination $\lambda$ on the surface receives a \emph{length}, which can be computed by measuring the (weighted) lengths of weighted curves converging to $\lambda$.
It is known \cite{bridgeman} that the lengths of the pleating laminations of a quasifuchsian group are bounded by a constant depending only on the underlying surface (here the punctured torus). By Thurston's double limit theorem (see Theorems 4.1 and 6.3 of \cite{thurstonlimitedouble}), the space of discrete, type-preserving representations in which two fixed measured laminations have length bounded by a given constant is \emph{compact}. Therefore, up to taking a subsequence, the groups $\Gamma^T$ corresponding to $w^T$ converge algebraically to a certain $\Gamma$. By \cite{minsky}, $\mathbb{H}^3/\Gamma$ is homeomorphic to $S\times \mathbb{R}$ and must have an infinite end (otherwise, $\Gamma$ would be quasifuchsian and the volumes would stay bounded). 

\begin{proposition}
The $\lambda^-$-end of $\mathbb{H}^3/\Gamma$ is finite, with pleating lamination $\lambda^-$. The $\lambda^+$-end is infinite, with projective ending lamination $[\lambda^+]$. \end{proposition}
\begin{proof}
By \cite{bonahon} (Theorem D), the space of type-preserving representations of the abstract group $\Gamma$ is smoothly (in fact, holomorphically) parametrized by the data $(\tau, \omega)$ of a point $\tau$ of Teichm\"uller space $\mathcal{T}$, and a \emph{transverse $\mathbb{R}/2\pi \mathbb{Z}$-valued cocycle} relative to a fixed topological lamination $\mu$ (such cocycles include pleating measures as special cases). Taking for $\mu$ the support of $\lambda^-$, we see that the moduli of the $\lambda^-$-boundaries of the convex cores of the $\Gamma^T$ must converge in $\mathcal{T}$. Therefore, $\mathbb{H}^3/\Gamma$ contains a locally convex pleated surface $H$ with pleating $\lambda^-$, which must be a boundary of the convex core ($\partial H$ contains all parabolic fixed points).

The $\lambda^+$-end of $\mathbb{H}^3/\Gamma$ must therefore be infinite. 
The parabolic fixed points of the limit group $\Gamma$ determine a version of Figure \ref{frise} (a Euclidean cusp link), and therefore a family of \emph{non-negative} angle assignments for the tetrahedra $\Delta_i$. By algebraic convergence, the $w^T$ converge to $w$ in the product topology. For any $i\leq j$, the total volume of $\Delta_{i-1}, \Delta_i, \dots, \Delta_j, \Delta_{j+1}$ is maximal with respect to $w_i, \dots,w_j$: in particular, Propositions \ref{tricot} and \ref{hingefragile} are still true. The techniques of \cite{mapomme} (Lemma 16 and Section 9 there) show that all $\Delta_i$ have positive angles, and Proposition \ref{tiltlemma}, hence also Lemma \ref{minkonwex}, still hold: $\{\Delta_i\}_{i\in\mathbb{Z}}$ is the geometrically canonical decomposition of $\mathbb{H}^3/\Gamma$. In particular, the family of all edges of all tetrahedra $\{\Delta_i\}_{i\geq 0}$ forms a sequence of laminations which exits $\mathbb{H}^3/\Gamma$: therefore $[\lambda^+]$ is the end invariant.
\end{proof}
The case of two infinite ends is already treated in \cite{akiyoshi}. Theorem \ref{ephhyper} is proved.

\subsection{The Pleating Lamination Conjecture for punctured-torus groups}

\begin{proposition} The group $\Gamma$ constructed at the end of Section \ref{sectionhyperbolicvolume} is \emph{continuously} parametrized by $(\lambda^+,\lambda^-)$. \llabel{groupcontinuous} \end{proposition}
\begin{proof}
Our first observation is that if $\beta^+/\alpha^+$ is rational, the initial choice of infinite prefix/suffix in Subsection \ref{onerationalpleating} does not change the resulting group $\Gamma$: it just induces different triangulations of the toric piece of Figure \ref{toricpieces}, whose deformation space is still the same.

Define the open set $$\mathcal{U}:=\mathbb{R}^2\smallsetminus \bigcup_{m,n\in \mathbb{Z}}[\pi,+\infty)\cdot\{(m,n)\}$$ (note that $0\notin \mathcal{U}$). An admissible pleating lamination $\lambda$ can be identified with an element $\pm(\alpha,\beta)$ of $\mathcal{U}/\pm$. Suppose $(\alpha^+_n, \beta^+_n) \rightarrow (\alpha^+, \beta^+)$ and $(\alpha^-_n, \beta^-_n) \rightarrow (\alpha^-, \beta^-)$ in $\mathcal{U}$, and define the oriented line $\Lambda_n$ from $\beta^-_n/\alpha^-_n$ to $\beta^+_n/\alpha^+_n$ across the Farey diagram. Also define the associated functions $\phi^{\pm,n}$ as in (\ref{phiphi}), domains $W^n$ as in Definition \ref{defineW}, and solutions $(w^n_i)_{i\in\mathbb{Z}}$ to the volume maximization problem over $W^n$. If $\beta^+/\alpha^+$ is rational, we may assume (up to restricting to two subsequences) that the $(\alpha^+_n, \beta^+_n)$ converge to $(\alpha^+,\beta^+)$ in the clockwise direction for the natural orientation of $\mathbb{P}^1\mathbb{R}$. We also make a similar assumption for $\beta^-/\alpha^-$.

A priori, the sequences $(w^n_i)_{i\in \mathbb{Z}}$ are defined only up to a shift of the index $i$. However, we can choose these shifts in a consistent way:
there exists a Farey edge $e$ which is crossed by all the lines $\Lambda_n$ for $n$ large enough, so we decide that $w^n_0$ always lives on $e$ (namely, $e^n_0=e$). By compactness of $[0,\pi]^{\mathbb{Z}}$, some subsequence of $(w^n)_{n\in \mathbb{N}}$ converges to some $w^*$. It is enough to show that $w^*=w$: indeed, the group $\Gamma$ is completely determined by the shapes of a finite number of tetrahedra $\Delta_i$.

The main observation is that the words $\Omega^n \in \{R,L\}^{\mathbb{Z}}$ converge pointwise to $\Omega$, and $\phi^{\pm,n}\rightarrow \phi^{\pm}$ (pointwise in $\mathbb{R}^{\mathbb{Z}}$), by definition (\ref{phiphi}). Therefore, $w^*$ belongs to the space $W$, hence the volume inequality $\mathcal{V}(w^*)\leq \mathcal{V}(w)$. Since $\max_W \mathcal{V}$ is achieved at a unique point (the volume is a strictly concave function), it is enough to prove the reverse inequality. We proceed by contradiction.

Suppose $\mathcal{V}(w^*) < \mathcal{V}(w)$. Pick $\varepsilon>0$: there exist integers $m<0<M$ such that the tetrahedra $\{\Delta_i\}_{m<i<M}$ defined by $w$ have total volume at least $\mathcal{V}(w)-\varepsilon$. If we can extend $(w_i)_{m\leq i\leq M}$ to a sequence $(v_i)_{i\in\mathbb{Z}}$ of $W^n$ for some large $n$, we will obtain a contradiction (assuming $\varepsilon$ small enough).

By Corollary \ref{hiatus2}, we can assume $\frac{w_m}{\phi^-_m}>\frac{w_{m+1}}{\phi^-_{m+1}}$ and $\frac{w_M}{\phi^-_M}>\frac{w_{M-1}}{\phi^-_{M-1}}$. Since $w$ satisfies (\ref{doubleve}) (strong inequalities), for $n$ large enough, the restricted sequence $(w_m,\dots,w_M)$ satisfies the corresponding inequalities defining $W^n$, by convergence of the $\phi^{\pm, n}$. Pick such a large $n$ and define
$$v_i:=\left \{ \begin{array}{ll} \phi^{-,n}_i \frac{w_m}{\phi^{-,n}_m} & \text{ if } i\leq m~; \\ w_i & \text{ if } m\leq i \leq M ~; \\ \phi^{+,n}_i \frac{w_M}{\phi^{+,n}_M} & \text{ if } M\leq i~. \end{array} \right .$$ It is a straightforward exercise to check that $(v_i)_{i\in \mathbb{Z}}$ belongs to $W^n$. \end{proof}

\begin{corollary} \emph{(C. Series)}\llabel{plth}
A quasifuchsian punctured-torus group $\Gamma$ is determined up to conjugacy in $\text{Isom}(\mathbb{H}^3)$ by its pleating measures $\lambda^{\pm}$.
\end{corollary}
\begin{proof}
We will use the well-known fact that the space $\mathcal{QF}$ of quasifuchsian, non-fuchsian (punctured-torus) groups is a connected real manifold of dimension $4$. Recall the open set $\mathcal{U}$ from the proof of Proposition \ref{groupcontinuous}, and consider the map
$$f:\mathcal{U} \times \mathcal{U} \longrightarrow \mathcal{QF}$$ defined by the construction of the group $\Gamma$ (end of Section \ref{sectionhyperbolicvolume}). We know that $f$ is well-defined and injective (Theorem \ref{pureexistence}), and continuous (Proposition \ref{groupcontinuous}). Since $\mathcal{U}^2$ has dimension $4$, the theorem of domain invariance states that the image $\text{Im}(f)$ is open. It remains to show that $\text{Im}(f)$ is closed.

Consider pairs $(\lambda^+_n, \lambda^-_n)$ such that the corresponding groups $\Gamma_n = f(\lambda^+_n, \lambda^-_n)$ converge to some $\Gamma$ in $\mathcal{QF}$. The function which to a group associates its pleatings is continuous (see \cite{keenseries}), so the $\lambda^{\pm}_n$ converge to some $\lambda^{\pm}$ in $\mathcal{U}$. Proposition \ref{groupcontinuous} then implies that $\Gamma=f(\lambda^+,\lambda^-)$. \end{proof}
Theorem \ref{ephgeneral} now follows from Theorem \ref{ephhyper}.

\begin{flushright}
USC Mathematics (KAP) \\ 3620 South Vermont Avenue \\ Los Angeles, CA 90089 (USA) \\
\texttt{Francois.Gueritaud@normalesup.org}
\end{flushright}

\begin{thebibliography}{MaR}

\bibitem[Ak]{akiyoshi} Hirotaka Akiyoshi, \emph{On the Ford domains of once-punctured torus groups}, in \emph{Hyperbolic spaces and related topics}, RIMS, Kyoto, Kokyuroku {\bf 1104} (1999), 109-121. 

\bibitem[AS]{comparing} Hirotaka Akiyoshi, Makoto Sakuma,
\emph{Comparing two convex hull constructions of cusped hyperbolic
manifolds}, Proceedings of the Workshop ``Kleinian groups and
hyperbolic 3-manifolds'' (Warwick 2002), Lond. Math. Soc. Lecture
Notes {\bf 299} (2003), 209--246.

\bibitem[ASWY1]{aswy1} H. Akiyoshi, M. Sakuma, M. Wada, Y. Yamashita, \emph{Jorgensen's picture of quasifuchsian punctured torus groups}, Proceedings of the Workshop ''Kleinian groups and hyperbolic 3-manifolds'' (Warwick 2002), Lond. Math. Soc. Lecture Notes {\bf 299} (2003), 247--273.

\bibitem[ASWY2]{aswy2} H. Akiyoshi, M. Sakuma, M. Wada, Y. Yamashita, \emph{Ford domains of punctured torus groups and two-bridge knot groups}, in \emph{Knot Theory}, Proceedings of the workshop held in Toronto dedicated to 70th birthday of Prof. K. Murasugi, 1999. 

\bibitem[ASWY3]{aswy3} Hirotaka Akiyoshi, Makoto Sakuma, Masaaki Wada,
Yasushi Yamashita, \emph{Punctured torus groups and 2-bridge knot groups} , Preprint.


\bibitem[Bon]{bonahon} Francis Bonahon, \emph{Shearing hyperbolic surfaces, bending pleated surfaces and Thurston's symplectic form}, Ann. Fac. Sci. Toulouse Math. {\bf 5} (1996),  233--297.

\bibitem[Bow]{bowditch} Brian H. Bowditch, \emph{Markoff triples and quasifuchsian groups}, Proc. Lond. Math. Soc. Vol. {\bf 77} (1998) 697--736.

\bibitem[Br]{bridgeman} Martin Bridgeman, \emph{Average bending of convex pleated planes in hyperbolic three-space}, Invent. Math. {\bf 132}, no. 3 (1998), 381--391.

\bibitem[CEG]{lms} R. D. Canary, D.B.A. Epstein, R. Green, \emph{Notes on notes of Thurston}, in
\emph{Analytical and Geometric Aspects of Hyperbolic Space} (Epstein
ed.), Lond. Math. Soc. Lecture notes {\bf 111}, Cambridge University
Press, 1987.

\bibitem[CH]{chanhodgson} Ken Chan, \emph{Constructing hyperbolic 3-manifolds}, Undergraduate thesis with Craig Hodgson, University of Melbourne, 2002.

\bibitem[EP]{epsteinpenner} David B.A. Epstein, Robert C. Penner, \emph{Euclidean decompositions of noncompact hyperbolic manifolds}, J. Diff. Geom. {\bf 27} (1988), 67--80.

\bibitem[G1]{abeilles} Fran\c{c}ois Gu\'eritaud, \emph{Formal Markoff maps are positive}, in preparation.

\bibitem[G2]{these} Franois Gu\'{e}ritaud, PhD thesis, in preparation.

\bibitem[GF]{mapomme} Fran\c{c}ois Gu\'eritaud, with an Appendix by David Futer,
\emph{On canonical triangulations of once-punctured torus bundles
and two-bridge link complements}, arXiv:math.GT/0406242, 2004.

\bibitem[J{\o}]{jorgensen} Troels J{\o}rgensen, \emph{On pairs of punctured tori}, (unfinished manuscript), in \emph{Kleinian groups and hyperbolic 3-manifolds} (Y. Komori, V. Markovic, C. Series, Editors), Lond. Math. Soc. Lecture notes {\bf 299} (2003), 183--207.

\bibitem[KS]{keenseries} Linda Keen and Caroline Series, \emph{Continuity of convex hull boundaries}, Pacific J. Math. {\bf 168}, no. 1 (1995), 183--206.

\bibitem[La]{lackenby} Marc Lackenby, \emph{The canonical decomposition of  once-punctured torus bundles}, Comment. Math. Helv. {\bf 78} (2003), 363--384.

\bibitem[Mil]{milnor} John Milnor, \emph{Hyperbolic geometry: the first 150
years}, Bull. Amer. Math. Soc. {\bf 6} (1982), no. 1, 9--24.

\bibitem[Min]{minsky} Yair Minsky, \emph{The classification of punctured-torus groups}, Annals of Math. {\bf 149} (1999), 559--626.


\bibitem[Ri]{rivin} Igor Rivin,
\emph{Euclidean structures on simplicial surfaces and hyperbolic volume},
Ann. of Math. {\bf 139} (1994) 553--580.

\bibitem[Se]{series} Caroline Series, \emph{Thurston's bending measure
conjecture for once punctured torus groups}, arXiv:math.GT/0406056,
2004.

\bibitem[Th]{thurstonlimitedouble} William P. Thurston, \emph{Hyperbolic Structures on 3-manifolds, II: Surface groups and 3-manifolds which fiber over the circle}, arXiv:math.GT/9801045.

\bibitem[Wa]{opti} Masaaki Wada, \emph{Opti}, a computer program for the study of punctured-torus groups: \texttt{http://vivaldi.ics.nara-wu.ac.jp/\~{}wada/OPTi/}

\end{thebibliography}
\end{document}